\documentclass{amsart}
\usepackage{amsmath,amssymb,amsxtra,amsthm,amscd}

\usepackage{lmodern}
\usepackage[T1]{fontenc}
\usepackage[all,cmtip]{xy}
\usepackage{rotating}
\usepackage{microtype}
\usepackage[color=orange,bordercolor=orange,linecolor=orange]{todonotes}
\usepackage{adjustbox}
\usepackage{xcolor}
\newif\ifShowLabels
\ShowLabelstrue
\newcommand{\TeXref}[1]{
\marginpar{\scriptsize \texttt{#1}}}
\ShowLabelsfalse

\DeclareMathOperator{\B}{\mathbf{B}}

         \newcommand{\bhlf}{\sideset{^{b\!}}{^{\textit{lf}}}\h}

\DeclareMathOperator{\C}{\mathbf{C}}

\DeclareMathOperator{\EGamma}{\boldsymbol{E}\boldsymbol{\Gamma}}

\DeclareMathOperator{\Free}{\mathbf{Free}}

\DeclareMathOperator{\Fun}{Fun}
\DeclareMathOperator{\G}{\mathbf{G}}
\DeclareMathOperator{\iG}{\mathit{i}\mathbf{G}}

         \newcommand{\Gnc}{G^{-\infty}}

         \newcommand{\GncX}{{G}^{-\infty}_{X}}

\DeclareMathOperator{\h}{\mathit{h}}
\DeclareMathOperator{\hlf}{\h^{\textit{lf}}}

\DeclareMathOperator{\Hom}{Hom}

\DeclareMathOperator{\id}{id}

\DeclareMathOperator{\im}{im}

\DeclareMathOperator{\J}{\mathit{J}}
         \newcommand{\Jh}{\J^{h}}
\DeclareMathOperator{\K}{\mathit{K}}
         \newcommand{\Knc}{\K^{-\infty}}

\DeclareMathOperator{\Mod}{\mathbf{Mod}}

\DeclareMathOperator{\point}{point}

\DeclareMathOperator{\Spt}{Spt}

         \newcommand{\subdot}{\boldsymbol{\cdot}}

\DeclareMathOperator{\U}{\mathbf{U}}

\DeclareMathOperator*{\one}{1}
\newcommand{\onehatplace}[1]
{ \one^{\substack{#1 \\ \frown}} }

\DeclareMathOperator*{\bones}{\times}
\newcommand{\undertimes}[1]
{ \bones_{#1} }

\DeclareMathOperator*{\bowl}{\cup}
\newcommand{\undercup}[1]
{ \bowl_{#1} }

\DeclareMathOperator*{\arch}{\cap}
\newcommand{\undercap}[1]
{ \arch_{#1} }

\newcommand{\pull}
{\!\!\! -\!\!\! -\!\!\! -\!\!\!}

\DeclareMathOperator*{\holimprep}{holim}                       
\newcommand{\holim}[1]%
{\displaystyle\holimprep_{\substack{\leftarrow \pull - \\ #1}} \, }

\DeclareMathOperator*{\hocolimprep}{hocolim}                   
\newcommand{\hocolim}[1]%
{\displaystyle\hocolimprep_{\substack{- \pull \rightarrow \\ #1}} \, }

\DeclareMathOperator*{\plainlim}{lim}                           
\newcommand{\contralim}[1]%
{\displaystyle\plainlim_{\substack{\leftarrow \pull - \\ #1}} \, }

\DeclareMathOperator*{\plaincolim}{colim}                       
\newcommand{\colim}[1]%
{\displaystyle\plaincolim_{\substack{- \pull \rightarrow \\ #1}} \, }

\DeclareMathOperator*{\laxlimplain}{laxlim}                     
\newcommand{\laxlim}[1]%
{\displaystyle\laxlimplain_{\substack{\leftarrow \pull - \\ #1}} \, }

\providecommand{\bysame}{\makebox[3em]{\hrulefill}\thinspace}

\theoremstyle{plain}
\newtheorem{Thm}{Theorem}[section]

\newtheorem{Cor}[Thm]{Corollary}

\newtheorem{Lem}[Thm]{Lemma}
\newtheorem{Prop}[Thm]{Proposition}

\theoremstyle{definition}
\newtheorem{Def}[Thm]{Definition}

\newtheorem{Ex}[Thm]{Example}

\newtheorem{Rem}[Thm]{Remark}

\theoremstyle{remark}
\newtheorem{Not}[Thm]{Notation}

\newtheoremstyle{freestylethm}{6pt}{6pt}{\itshape}{}%
                {\bfseries}{}{.5em}{\thmnote{#3}}
\theoremstyle{freestylethm}

\newcommand{\SecRef}[2]{\section{#1}\label{S:#2}%
\ifShowLabels \TeXref{{S:#2}} \fi}
\newcommand{\SSecRef}[2]{\subsection{#1}\label{SS:#2}%
\ifShowLabels \TeXref{{SS:#2}} \fi}

\newcommand{\refS}[1]{\textup{\ref{S:#1}}}

\newcommand{\refT}[1]{\textup{\ref{T:#1}}}
\newcommand{\refL}[1]{\textup{\ref{L:#1}}}
\newcommand{\refD}[1]{\textup{\ref{D:#1}}}

\newcommand{\refP}[1]{\textup{\ref{P:#1}}}
\newcommand{\refR}[1]{\textup{\ref{R:#1}}}
\newcommand{\refN}[1]{\textup{\ref{N:#1}}}

\newenvironment{ThmRef}[1]%
{ \begin{Thm} \label{T:#1}
\ifShowLabels \TeXref{T:#1} \fi }%
{ \end{Thm} }
\newenvironment{DefRef}[1]%
{ \begin{Def} \label{D:#1}
\ifShowLabels \TeXref{D:#1} \fi }%
{ \end{Def} }
\newenvironment{LemRef}[1]%
{ \begin{Lem} \label{L:#1}
\ifShowLabels \TeXref{L:#1} \fi }%
{ \end{Lem} }
{ \begin{Cor} \label{C:#1}
\ifShowLabels \TeXref{C:#1} \fi }%
{ \end{Cor} }
\newenvironment{RemRef}[1]%
{ \begin{Rem} \label{R:#1}
\ifShowLabels \TeXref{R:#1} \fi }%
{ \end{Rem} }
\newenvironment{PropRef}[1]%
{ \begin{Prop} \label{P:#1}
\ifShowLabels \TeXref{P:#1} \fi }%
{ \end{Prop} }
\newenvironment{ExRef}[1]%
{ \begin{Ex} \label{E:#1}
\ifShowLabels \TeXref{E:#1} \fi  }%
{ \end{Ex} }
\newenvironment{NotRef}[1]%
{ \begin{Not} \label{N:#1}
\ifShowLabels \TeXref{N:#1} \fi }%
{ \end{Not} }

\newenvironment{ThmRefName}[2]%
{ \begin{Thm} [#2]\label{T:#1}
\ifShowLabels \TeXref{T:#1} \fi }%
{ \end{Thm} }
\newenvironment{DefRefName}[2]%
{ \begin{Def} [#2]\label{D:#1}
\ifShowLabels \TeXref{D:#1} \fi }%
{ \end{Def} }
{ \begin{Lem} [#2]\label{L:#1}
\ifShowLabels \TeXref{L:#1} \fi }%
{ \end{Lem} }
{ \begin{Cor} [#2]\label{C:#1}
\ifShowLabels \TeXref{C:#1} \fi }%
{ \end{Cor} }
{ \begin{Rem} [#2]\label{R:#1}
\ifShowLabels \TeXref{R:#1} \fi }%
{ \end{Rem} }
{ \begin{Prop} [#2]\label{P:#1}
\ifShowLabels \TeXref{P:#1} \fi }%
{ \end{Prop} }
\newenvironment{ExRefName}[2]%
{ \begin{Ex} [#2]\label{E:#1}
\ifShowLabels \TeXref{E:#1} \fi }%
{ \end{Ex} }

\let\oldtocsection=\tocsection
\let\oldtocsubsection=\tocsubsection
\renewcommand{\tocsection}[2]{\hspace{0em}\oldtocsection{#1}{#2}}
\renewcommand{\tocsubsection}[2]{\hspace{2em}\oldtocsubsection{#1}{#2}}

\begin{document}

\title{Excision in equivariant fibred \textit{G}-theory}
\author[Gunnar Carlsson]{Gunnar Carlsson}
\address{Department of Mathematics\\ Stanford University\\ Stanford\\ CA 94305}
\email{gunnar@math.stanford.edu}
\author[Boris Goldfarb]{Boris Goldfarb}
\address{Department of Mathematics and Statistics\\ SUNY\\ Albany\\ NY 12222}
\email{goldfarb@math.albany.edu}
\date{\today}

\begin{abstract}
This paper provides a generalization of excision theorems in controlled algebra in the context of equivariant $G$-theory with fibred control and families of bounded actions.  It also states and proves several characteristic features of this theory such as existence of the fibred assembly and the fibrewise trivialization.
\end{abstract}

\maketitle

\tableofcontents

\SecRef{Introduction}{INTRO}

The bounded $K$-theory construction due to Pedersen and Weibel \cite{ePcW:85} has been shown to be extremely useful in the analysis of versions of the Novikov conjecture (\cite{gC:95,gCbG:04,gCeP:93,gCeP:97,dRrTgY:14}).  This conjecture asserts the split injectivity of a natural transformation called the {\em assembly}.  The present paper is the culmination of a series of papers of the authors \cite{gCbG:00,gCbG:18a,gCbG:18b}  that extend the techniques sufficient to address the much more difficult Borel conjecture, which is very closely related to the question of whether the assembly map is an isomorphism.  What we have found is that substantial extensions are necessary.  In order to explain these extensions, we first recall some of  the properties of the original construction of Pedersen and Weibel.  
\medskip
\begin{itemize}
\item{The construction begins with a metric space $X$ and a commutative ring $R$ of coefficients and constructs a bounded $K$-theory spectrum $K(X,R)$.}  {The spectrum is constructed by considering potentially infinitely generated based $R$-modules equipped with a reference map $\varphi$ from the chosen basis $B$ to the metric space $X$, and considering only morphisms satisfying a condition of control, defined by a single parameter $P$. The requirement for a homomorphism $f$  from $(M,B_M,\varphi_M) $ to $(N, B_N, \varphi _N)$ to be bounded with parameter $P$ is that for any $b \in B_M$, $f(b)$ lies in the span of the basis elements in $B_N$ which lie within a distance $P$ of $\varphi _M(b)$. We will call such modules {\em geometric $R$-modules over $X$}.} \medskip
\item{The construction is functorial for maps $f \colon X \to Y$ of metric spaces that are {\em unifomly expansive}, which are maps such that there is a function $c \colon \Bbb{R} \rightarrow \Bbb{R}$ so that if $x,x^{\prime} \in X$ are any two points with $d(x,x^{\prime}) \leq t$, then $d(f(x),f(x^{\prime})) \leq c(t)$.  } \medskip
\item{The construction is {\em coarse invariant}.  For example, if we have an isometric  inclusion $X\hookrightarrow Y$ of metric spaces, and every element of $Y$ is within a fixed distance of an element of $X$, then the map $K(X,R) \rightarrow K(Y,R)$ is an equivalence of spectra.  In particular, it is possible to show that for compact  $K(\Gamma,1)$-manifolds $X$, the inclusion of a single orbit $\Gamma \cdot x$ in the universal cover of $X$ is an equivalence.} \medskip
\item{For $X$ a Riemannian manifold, the construction is closely related to the {\em locally finite} or {\em Borel homology} of $X$ with coefficients in the algebraic $K$-theory spectrum $K(R)$. Specifically, there is a bounded assembly map $$\hlf (X,K(R)) \longrightarrow K(X,R)$$ which is an equivalence for a large class of manifolds.   } \medskip
\item{For free and proper actions by isometries of a discrete group on a Riemannian manifold $X$, there is an equivariant construction $K^{\Gamma,0}(X,R) $ explained in section 3 which is equivalent to the original construction $K(X,R)$ non-equivariantly, and whose fixed point spectrum is the ordinary homology of $X/\Gamma _+$ with coefficients in $K(R)$.  }
\end{itemize} 

\medskip

There are various directions in which one can generalize the Pedersen-Weibel framework. The following are the extensions and issues that need to be addressed, with some comments and motivation. 
\medskip

\begin{itemize}
\item{In algebraic $K$-theory, one usually considers categories of free or projective finitely generated modules, however in many situations one can use other categories of modules such as finitely generated or finitely presented modules.  When studying Noetherian commutative rings, one can consider the category of finitely generated modules and perform an algebraic $K$-theoretic construction on it.  The result which is often easier to compute is referred to as {\em $G$-theory}.  The authors have performed analogues of this kind of construction in the context of bounded $K$-theory and group rings in \cite{gCbG:00,gCbG:15}. } \medskip
\item{The issue that naturally comes up is construction of equivariant theories beyond the case of isometric actions.  We are led to study actions of a group $\Gamma$ on a metric space through uniformly expansive self-maps of metric spaces.  These actions can no longer be assumed to be free, and each self-map is no longer necessarily an injection.}\medskip
\item{Next, there is the related new requirement of more general notions of control.  Specifically, given two metric spaces $X$ and $Y$, we will need to consider the notion of {\em fibrewise control} for geometric $R$-modules over $X \times Y$. Given two  geometric modules $M = (M,B_M, \varphi _M) $ and $N = (N, B_N, \varphi _N) $  over  $X \times Y$, equipped with the sum metric $D = d_X + d_Y$,  we say a homomorphism $f \colon M \rightarrow N$ is {\em fibrewise controlled} if it satisfies two conditions.
\begin{enumerate}
\item{There exists a number $P>0$ so that for any basis element $b \in B_M$, $f(b)$ lies in the span of the those basis elements $b^{\prime}$  in $B_N$ for which $\pi _X(\varphi(b^{\prime}))$ lies in a ball of radius $P$ of $\pi _X(\varphi_M(b))$ in $X$.   }
\item{For any bounded subset $U$ in $X$, there is a number $P(U) > 0$ so that for any $b \in (\pi _X \circ \varphi _M)^{-1}(U)$, $f(b)$ lies in the span of the set of elements $b^{\prime}$  of $B_N$ for which $d(\pi _Y\circ\varphi _M(b), \pi _Y \circ \varphi _N(b^{\prime})) \leq P(U)$.   }
\end{enumerate}
Intuitively, this notion of control provides maps which are controlled in the ``base direction" in the bundle $X \times Y \rightarrow X$, and which are controlled in each of the "fibers", but not necessarily requiring that there exist a uniform bound over all fibers.  We will also require $G$-theoretic and equivariant versions of such  notions of control.    } \medskip
\item{Fibrewise control is the analogue for $K$- and $G$-theory of the notion of ``parametrized homotopy theory" or ``homotopy theory over a base", where $X$ is the base and $Y$ is the fiber.  In the equivariant case, where $X$ is equal to a group $\Gamma$ regarded as a metric space using the word length metric.  The fixed points of the $\Gamma$ action on $K_{\Gamma}(Y)$  is analogous to the study of the bundles over a classifying space $B\Gamma$ obtained from $\Gamma$-spaces $Y$ by the construction $Y \rightarrow E\Gamma {\times}_{\Gamma} Y$. To realize our results, we will need excision properties holding for coverings of $Y$.   } \medskip
\item{Finally, we will need results that imply that $K_{\Gamma}(Y) $ is canonically equivalent to $K_{\Gamma}(Y^0)$, where $Y^0$ is defined to be $Y$ equipped with the trivial action.    We say that an action of a group $\Gamma$ on a metric space is {\em bounded} if for every $\gamma \in \Gamma$, there exists an $R_{\gamma} > 0 $ such that $d(x,\gamma (x)) \leq R_{\gamma}$ for all $x \in X$.  We will be able to prove that for any bounded action of $\Gamma$ on a metric space $Y$, $K_{\Gamma}(Y)$ is canonically equivalent to $K_{\Gamma}(Y^0)$.  

That such a result should be possible is suggested by the analogue that occurs in ordinary parametrized homotopy theory, where $E\Gamma {\times}_{\Gamma} BH$ is canonically equivalent to $B\Gamma \times BH $ whenever the $\Gamma$-action is specified by a homomorphism to $H$ followed by the homomorphism $H \rightarrow \mathit{Inn}(H)$, where $\mathit{Inn}(H)$ is the inner automorphism group of $H$.  The reason this is an apt analogy is that in the case where we are considering a bounded action of $\Gamma$ on a metric space $Y$, the action can be considered as conjugation by bounded automorphisms of $Y$, all of which are contained in the category defining $K(Y)$.     }
\end{itemize} 

In this paper we prove excision results that incorporate all generalizations simultaneously. Because we will only need the excision results where the action on the space $Y$ is bounded in the sense defined above, we will only prove them in that situation.   That is, we prove excision theorems (Theorems \refT{Exc2} and \refT{PPPOI8fin}) for equivariant $G$-theory with fibred control of bounded $\Gamma$-spaces $Y$.  Additionally we obtain the results suggested above as part of equivariant fibred $G$-theory.

\SecRef{Homotopy fixed points in categories with action}{HFP}

Given an action of a group ${\Gamma}$ on a space $X$, one has the subspace of fixed points $X^{\Gamma}$.  This subspace often has geometric  significance for the study of $X$ and ${\Gamma}$.  A different powerful idea in topology is to model an interesting space or spectrum as the fixed point space or spectrum $X^{\Gamma}$ for a specifically designed $X$ with an action by a related group ${\Gamma}$.  In either case, there is always the homotopy fixed point spectrum $X^{h{\Gamma}}$ which is easier to understand than $X^{\Gamma}$ and the canonical reference map $\rho \colon X^{{\Gamma}} \to X^{h{\Gamma}}$.

Now suppose we have a group action on a category.  This automatically produces an action on the nerve and therefore a space.  Suppose the category is then fed into a machine such as the algebraic $K$-theory, and we are interested in the fixed points of the $K$-theory.  Therefore we want to look at the homotopy fixed points.  In many important cases it is possible to construct a spatial or categorical description of what we get.  Thomason defined the lax limit category whose $K$-theory turns out to be exactly the homotopy fixed points of the old action.

\begin{DefRef}{ThomaPrep}
	Let $\EGamma$ be the category with the object set ${\Gamma}$ and the
unique morphism $\mu \colon \gamma_1 \to \gamma_2$ for any pair
$\gamma_1$, $\gamma_2 \in {\Gamma}$. There is a left ${\Gamma}$-action
on $\EGamma$ induced by the left multiplication in ${\Gamma}$.
If $\mathcal{C}$ is a category with a left ${\Gamma}$-action,
then the category of functors
$\Fun(\EGamma,\mathcal{C})$ is another
category with the ${\Gamma}$-action given on objects by
the formulas $\gamma(F)(\gamma')=\gamma F (\gamma^{-1} \gamma')$
and $\gamma(F)(\mu)=\gamma F (\gamma^{-1} \mu)$.
It is
nonequivariantly equivalent to $\mathcal{C}$. 
\end{DefRef}

The category $\Fun(\EGamma,\mathcal{C})$ is an interesting and useful object in its own right.  There are several manifestations of this, for example in the work of Mona Merling and coauthors \cite{bGpMmM:17,cMmM:16,mM:16} or the work of these authors \cite{gC:95,gCbG:03,gCbG:04,gCbG:13}.  While in both applications it is crucial to work with the category itself, in this paper we concentrate on approximating the fixed points in $\Fun(\EGamma,\mathcal{C})$.

The following construction has been used by Thomason \cite{rT:83}.  We will refer to it as the \textit{homotopy fixed points of a category}, following Merling \cite{mM:16}.

\begin{DefRefName}{Thoma}{Homotopy fixed points}
The fixed point subcategory $\Fun(\EGamma,\mathcal{C})^{{\Gamma}}$ of the category of functors $\Fun(\EGamma,\mathcal{C})$ consists of equivariant functors and
equivariant natural transformations.  We will denote it by $\mathcal{C}^{h{\Gamma}}$.

Explicitly, the objects of $\mathcal{C}^{h{\Gamma}}$ are the pairs
$(C,\psi)$ where $C$ is an object of $\mathcal{C}$ and $\psi$ is a function from
${\Gamma}$ to the morphisms of $\mathcal{C}$ with $\psi (\gamma) \in \Hom (C,\gamma C)$ that satisfies 
$\psi(e) = \id$ for the identity group element $e$, and satisfies the cocycle identity $\psi (\gamma_1 \gamma_2) =
\gamma_1 \psi(\gamma_2)  \psi (\gamma_1)$ for all pairs $\gamma_1$ and $\gamma_2$ in ${\Gamma}$.
These conditions imply that $\psi (\gamma)$ is always an
isomorphism. The set of morphisms $(C,\psi) \to
(C',\psi')$ consists of the morphisms $\phi \colon C \to C'$ in
$\mathcal{C}$ such that the squares
\[
\xymatrix{
C \ar[rr]^-{\psi (\gamma)} \ar[dd]_-{\phi} &&\gamma C \ar[dd]^-{\gamma \phi}\\
\\
C' \ar[rr]^-{\psi' (\gamma)} &&\gamma C'
}
\]
commute for all $\gamma \in {\Gamma}$.
\end{DefRefName}

\begin{RemRef}{UOFFV}
	As pointed out in \cite{mM:16}, the homotopy fixed points of a category are not necessarily identical with space level constructions.  It is for example not true in general that the nerve of the homotopy fixed point category  of a category is the same as the geometric homotopy fixed points of the nerve of the category. It is however true in the case where the category is a discrete $\Gamma$-groupoid.
\end{RemRef}

\begin{Ex} \label{groupring}  Let $\mathcal{C}$ denote a category, and equip it with the trivial  action by $\Gamma$.  Then the category $\mathcal{C}^{h\Gamma}$ is the category of representations of $\Gamma$ in $\mathcal{C}$. In particular, if $\mathcal{C}$ is the category of $R$-modules for a commutative ring $R$, then  $\mathcal{C}^{h\Gamma}$ may be identified with  the category of (left) $R[\Gamma]$-modules. 
\end{Ex} 

\begin{Ex} \label{fields}  Let $F \subseteq E$ denote a Galois field extension, with Galois group $G$.  We consider the skew group ring $\Lambda = E^t[G]$, and consider the category $\mathcal{C}_E$  whose objects are $F^t[G]$-modules and whose morphisms are the $E$-linear maps.  There is a $G$ action on $\mathcal{C}_E$, which is the identity on objects and which is defined by the group action on the morphisms.  In this case, $\mathcal{C}_E^{h \Gamma}$ is equivalent to the category of $F$-vector spaces 
\end{Ex} 

In Example \ref{groupring}, we saw that the group ring of  a group $\Gamma$ with coefficients in a commutative ring $R$ may be realized as the fixed point subcategory of the action of $\Gamma$ on $\Fun(\EGamma,\mathcal{C})$, where $\mathcal{C}$ denotes the category of all $R$-modules.  In many cases, however, it is important to understand the category of free and finitely generated left $R[\Gamma ]$-modules as a fixed point category.  This is the case in the papers \cite{gC:05} and \cite{gCbG:04}, for instance, where the injectivity of the assembly map is proved in a large family of cases. In the case of these two papers, this  is achieved by defining a subcategory of   $\Fun(\EGamma,\mathcal{C})$ by restricting the morphisms $\psi (\gamma)$.  The restriction in this case arises by the selection of a subcategory of the category of all $R$-modules based on the Pedersen-Weibel construction, which is endowed with a filtration and an action of the group $\Gamma$.  The restricted version of $\mathcal{C}_E^{h \Gamma}$ requires that all of the morphisms $\psi(\gamma)$ have filtration zero.  In order to attack the surjectivity problem for the assembly, we are led to the construction of more general forms of restriction of the maps $\psi (\gamma)$. This leads us to the concept of the {\em relative homotopy fixed points of a category}, which we now define.   

\begin{DefRefName}{RLL}{Relative homotopy fixed points}
The category $\mathcal{C}^{h{\Gamma}} (\mathcal{M})$ is defined using input data consisting of a category $\mathcal{C}$  equipped with an action by a group $\Gamma$ and a subcategory $\mathcal{M} \subset \mathcal{C}$ closed under the action of $\Gamma$.  It is the full subcategory of $\mathcal{C}^{h \Gamma}$ on objects $(C, \psi )$ with the additional condition that $\psi(\gamma)$ is in $\mathcal{M}$ for all elements $\gamma \in \Gamma$. 
\end{DefRefName}

\begin{ExRef}{KLO}
	Clearly, if $\mathcal{M}$ is the entire category $\mathcal{C}$, the relative homotopy fixed points are the genuine homotopy fixed points.
\end{ExRef}

\begin{ExRef}{NOV} In the case where $\mathcal{C}$ is a filtered category, we can consider the situation where $\mathcal{M}$ is the subcategory of the filtration zero morphisms.  This is the situation used in \cite{gC:05} and \cite{gCbG:04}.  
\end{ExRef}

We will exploit the relative homotopy fixed points in two applications. 
The first construction required in \cite{gC:95} by the first author allows us to model the $K$-theory of a group ring whenever the group has a finite classifying space. It is based on bounded $K$-theory of the group given a word metric with the isometric action on itself given by the left multiplication.  It turns out that the categorical homotopy fixed point construction equires a constraint.  We review that construction in section \refS{BKTKTGR}. 

A generalization of bounded $K$-theory has appeared in recent work on the Borel isomorphism conjecture \cite{gCbG:13,gCbG:15,gCbG:18a,gCbG:18b}. The fibred homotopy fixed points defined in \cite{gCbG:18a} are particularly intriguing because the standard tools in controlled algebra, based on Karoubi filtrations that generate long exact sequences, fail to compute them on a very basic level.  The details can be found in \cite[Example 5.2]{gCbG:18a}.  
Section \refS{FFLLIA} contains the definition of fibred homotopy fixed points and a precise statement of the desired facts that are unknown.  At the end of the section we relate this construction to the Borel conjecture. 
The resolution is provided in the last section \refS{TTITS} of this paper in terms of fibred homotopy fixed points in $G$-theory which we are able to compute.

\SecRef{Bounded \textit{K}-theory and the \textit{K}-theory of group rings}{BKTKTGR}

Bounded control is the simplest version of a ``control condition'' that can be imposed in various categories of modules, to which one can apply the algebraic $K$-theory construction. It was introduced in Pedersen \cite{eP:84} and Pedersen/Weibel \cite{ePcW:85} and has become crucial for $K$-theory computations in geometric topology.

Let $X$ be a metric space and let $R$ be an arbitrary associative ring with unity.  We will always assume that metric spaces are proper in the sense that closed bounded subsets are compact. 

\begin{DefRef}{PWcats}
The objects of the category of \textit{geometric $R$-modules over $X$} are locally finite functions $F$ from points of $X$ to the category of finitely generated free $R$-modules $\Free_{\mathit{fg}} (R)$.  Following Pedersen and Weibel, we will denote by $F_x$ the module assigned to the point $x$ of $X$ and denote the object itself by writing down the collection $\{ F_x \}$.
The \textit{local finiteness} condition requires precisely that for every bounded subset $S \subset X$ the restriction of $F$ to $S$ has finitely many nonzero modules as values.

Let $d$ be the distance function in $X$.  The morphisms $\phi \colon \{ F_x \} \to \{ G_x \}$ are collections of $R$-linear homomorphisms
$\phi_{x,x'} \colon F_x \rightarrow G_{x'}$,
for all $x$ and $x'$ in $X$, with the property
that $\phi_{x,x'}$ is the zero homomorphism whenever $d(x,x') > D$
for some fixed real number $D = D (\phi) \ge 0$.
One says that $\phi$ is \textit{bounded by} $D$.
The composition of two morphisms $\phi \colon \{ F_x \} \to \{ G_x \}$ and $\psi \colon \{ G_x \} \to \{ H_x \}$
is given by the formula
\begin{equation}
(\psi \circ \phi)_{x,x'} = \sum_{z \in X} \psi_{z,x'} \circ \phi_{x,z}.
\tag{$\ast$}
\end{equation}
This sum is finite because of the local finiteness property of $G$.

We will want to enlarge this category, and so we use instead an equivalent category $\mathcal{B} (X,R)$ that is better for this purpose.  

The objects are functors $F \colon \mathcal{P}(X) \to \Free (R)$ from the power set $\mathcal{P}(X)$ to the category of free modules, both viewed as posets ordered by split inclusions.  There are two additional requirements.  For every bounded subset $C$ of $X$ the value $F(C)$ has to belong to the subcategory of finitely generated modules $\Free_{\mathit{fg}} (R)$.  In the codomain, the values are required to satisfy the equality $F(S) = \bigoplus_{x \in S} F(x)$ for all $S \subset X$.  The morphisms in this reformulation are $R$-linear homomorphisms $\phi \colon F(X) \to G(X)$ such that the components
$\phi_{x,x'} \colon F(x) \to G(x')$ are zero whenever $d(x,x') > D$
for some $D$.
The composition of two morphisms $\phi \colon F \to G$ and $\psi \colon G \to H$ is the usual composition of $R$-linear homomorphisms; its components are the maps $(\psi \circ \phi)_{x,x'}$ in the formula ($\ast$) above.
\end{DefRef}

\begin{DefRef}{Bornol2}
A map $f \colon X \to Y$
between metric spaces is called \textit{uniformly expansive} if there is a function $\lambda \colon [0,\infty) \to [0,\infty)$
such that
$d_X (x_1, x_2) \le r$ implies
$d_Y (f(x_1), f(x_2)) \le \lambda (r)$.
A map $f$ is \textit{proper} if $f^{-1} (S)$ is a bounded subset of $X$ for
each bounded subset $S$ of $Y$.
We say $f$ is a \textit{coarse map} if it is eventually continuous and proper.
\end{DefRef}

Extensively used instances of coarse maps in geometry are quasi-isometries.

It is elementary to check that the geometric $R$-modules over $X$ is an additive category and that
coarse maps between metric spaces induce additive functors.
A coarse map $f$ is a \textit{coarse equivalence} if there is a coarse map $g \colon Y \to X$ such
that $f \circ g$ and $g \circ f$ are bounded maps.  It follows that an action of a group on a metric space by coarse equivalences induces an additive action on $\mathcal{B} (X,R)$.

We will treat the group $\Gamma$ equipped with a finite generating set ${\Omega}$ closed under taking inverses as a metric space.
The \textit{word-length metric} $d = d_{\Omega}$ is
induced from the condition that $d (\gamma, \gamma \omega) =
1$ whenever $\gamma \in \Gamma$ and $\omega \in \Omega$.
It is well-known that varying $\Omega$ only changes $\Gamma$ to a
quasi-isometric metric space.
The word-length metric makes $\Gamma$ a proper metric space with a free $\Gamma$-action by isometries via left
multiplication.

\medskip

\textit{An important observation.}
A free action of $\Gamma$ on $X$ by isometries always gives a free action on $\mathcal{C} = \mathcal{B} (X,R)$.  In contrast, $\Fun(\EGamma,\mathcal{C})$ with the induced group action does have the subcategory $\mathcal{C}^{h\Gamma}$ of equivariant functors. These homotopy fixed points, however, is not the correct notion for modeling the finitely generated free modules over $R[G]$ and the $K$-theory of $R[G]$.  
 
\begin{DefRef}{NMEWQS}
The category $\mathcal{B}^{\Gamma, 0} (X,R)^{\Gamma}$ are relative homotopy fixed points $\mathcal{C}^{h\Gamma} (\mathcal{M})$ with the following data:
$\mathcal{C}$ is the category of geometric modules $\mathcal{B} (X,R)$ and $\mathcal{M}$ consists of those morphisms in $\mathcal{C}$ that are bounded by $0$.
\end{DefRef}

The additive category $\mathcal{B}^{\Gamma, 0} (X,R)^{\Gamma}$ has the associated nonconnective $K$-theory spectrum 
$\Knc  (X,R)^{\Gamma}$ constructed as in \cite{ePcW:85}. 
There is now the following desired identification.

\begin{ThmRef}{JGun}
Suppose $\Gamma$ acts on $X$ freely, properly discontinuously by isometries so that the orbit space $X/\Gamma$ with the orbit metric is bounded. 
  	
It follows that 
$\Knc  (X,R)^{\Gamma}$ is weakly homotopy equivalent to the nonconnective spectrum $\Knc (R[\Gamma])$.  The stable homotopy groups of the nonconnective spectrum are the Quillen $K$-groups of $R[G]$ in nonnegative dimensions and the negative $K$-groups of Bass in negative dimensions.
\end{ThmRef}

\begin{proof}
The result follows from Corollary VI.8 in \cite{gC:95}.
\end{proof}

This geometric situation occurs, for example, when $\Gamma$ acts cocompactly, freely properly discontinuously on a contractible connected Riemannian manifold $X$ or when it acts on itself with a word metric via left multiplication.

\SecRef{Fibred homotopy fixed points in $K$-theory}{FFLLIA}

Let $\mathcal{A}$ be an additive category. Generalizing Definition \refD{PWcats}, one has the bounded category with coefficients in $\mathcal{A}$.

\begin{NotRef}{NOT}
Given a subset $S$ of a metric space and a number $k \ge 0$, $S[k]$ is used for the 	$k$-\textit{enlargement} of $S$ defined as the set of all points $x$ with $d(x, S) \le k$.
\end{NotRef}

Recall that $\mathcal{A}$ is a subcategory of its cocompletion ${\mathcal{A}}^{\ast}$ which is closed under colimits.  For example, a construction based on the presheaf category was given by Kelly in 6.23 of \cite{gK:82}.

\begin{DefRef}{OPQWS}
	$\mathcal{B} (X,\mathcal{A})$ has objects which are covariant functors $F \colon \mathcal{P}(X) \to {\mathcal{A}}^{\ast}$ from the power set $\mathcal{P}(X)$ to ${\mathcal{A}}^{\ast}$, both ordered by inclusion.  Just as in Definition \refD{PWcats}, there are several requirements: 
\begin{itemize}
\item $F(x)$ is an object of $\mathcal{A}$ for every point $x$ in $X$,
\item the resulting function $F \colon X \to \mathcal{A}$ is locally finite, so only finitely many values are non-zero when restricted to any compact subset of $X$,
\item for all subsets $S \subset X$, \[ F(S) = \bigoplus_{x \in S} F(x), \] 
\item the inclusion $F(S \subset X)$ is onto a direct summand for each subset $S$.
\end{itemize}
A morphism in $\mathcal{B} (X,\mathcal{A})$ is a morphism $\phi \colon F(X) \to G(X)$ in ${\mathcal{A}}^{\ast}$ with a number $D \ge 0$ such that $\phi$ restricted to $F(S)$ factors through $G(S[D])$ for all $S \subset X$. We say a morphism which admits such a number $D$ is {\em  $D$-controlled}
\end{DefRef}

This context, which produces a category isomorphic to $\mathcal{B} (X,R)$ when 
$\mathcal{A}$ is the category of free finitely generated $R$-modules, allows us  to iterate the bounded control construction as follows.

\begin{DefRefName}{Fibred}{Fibred control for geometric modules}
Given two metric spaces $X$ and $Y$ and any ring $R$, the category 
$\mathcal{B}_X (Y,R)$, or simply $\mathcal{B}_X (Y)$ when the choice of ring $R$ is clear, is the bounded category $\mathcal{B}(X, \mathcal{A})$ with $\mathcal{A} = \mathcal{B}(Y,R)$.
\end{DefRefName}

Among many options for relativizing homotopy fixed points in this setting, there is one of specific interest. 

Let $\mathcal{A} = \mathcal{B}(Y,R)$ as before and $\mathcal{A}' = \Mod (R)$ be the category of arbitrary $R$-modules.  There is a \textit{forget control} functor $t \colon \mathcal{B}(Y,R) \to \Mod (R)$ which only remembers that the objects are $R$-modules and the morphisms are $R$-linear homomorphisms.  From $t$ we may induce the functor $T \colon \mathcal{B}(X, \mathcal{A}) \to \mathcal{B}(X, \mathcal{A}')$.  

For this construction we assume that $\Gamma$ acts on $X$ by isometries and so, therefore, on $\mathcal{B}(X, \mathcal{A}')$. 
On the other hand, we allow the action of $\Gamma$ on $Y$ to be by coarse equivalences.  This can also be used to induce an action on $\mathcal{B}(X, \mathcal{A})$.

\begin{DefRefName}{RLL2}{Fibred homotopy fixed points in bounded $K$-theory}
These are relative homotopy fixed points with the following choice of ingredients: 
\begin{itemize}
	\item[--] the category $\mathcal{C}$ is $\mathcal{B}_X (Y,R)$, 
	\item[--] the subcategory $\mathcal{M}$ consists of all morphisms $\phi$ such that $T(\phi)$ is a controlled morphism bounded by $0$.
\end{itemize}
\end{DefRefName}

\begin{NotRef}{NMDKV}
When $X$ is the group $\Gamma$ itself with the left multiplication action and the word metric with respect to some choice of a finite set of generators, we will use the special notation $\mathcal{B}^{w\Gamma} (Y)$ for the fibred homotopy fixed points $\mathcal{C}^{h\Gamma} (\mathcal{M})$. 
\end{NotRef}

\SecRef{Summary of bounded \textit{G}-theory with fibred control}{KICVB}

A comprehensive exposition of bounded \textit{G}-theory with fibred control is available in \cite{gCbG:18b}.  This is a summary of that theory and a number of facts in the form we can refer to in the next section.

Throughout the rest of the paper, $R$ will be a Noetherian ring.

At the basic level, bounded $G$-theory with fibred control is an analogue of the algebraic $K$-theory of $\mathcal{B}_X (Y,R)$ locally modeled on finitely generated $R$-modules.  The result is an exact category $\B_X (Y)$ where the exact sequences are not necessarily split but which contains $\mathcal{B}_X (Y)$ as an exact subcategory.  

\begin{DefRef}{GEPCW2}
Given an $R$-module $F$, an $(X,Y)$\textit{-filtration} of $F$ is a functor
$ \phi_F \colon \mathcal{P}(X \times Y) \to  \mathcal{I}(F)$
from the power set of the product metric space to the partially ordered family of $R$-submodules of $F(X \times Y)$, both ordered by inclusion.
When there is no ambiguity, we find it convenient to denote the values $\phi_F (U)$ by $F(U)$.
We assume that $F$ is reduced
in the sense that $F(\emptyset)=0$.

The associated $X$-filtered $R$-module $F_X$ is given by $ F_X (S) = F(S \times Y)$.  Similarly, for each subset $S \subset X$, one has the $Y$-filtered $R$-module $F^S$ given by $F^S (T) = F(S \times T)$.
In particular, $F^X(T) = F(X \times T)$.
\end{DefRef}

We will use the following notation generalizing enlargements in a metric space.
Given a subset $U$ of $X \times Y$ and a function $k \colon X \to [0, + \infty )$, let
\[
U [k] = \{ (x,y) \in X \times Y \ \vert \ \textrm{there\ is} \ (x,y') \in U \ \textrm{with} \ d(y,y') \le k(x) \}.
\]
If in addition we are given a number $K \ge 0$ then
\[
U [K,k] = \{ (x,y) \in X \times Y \ \vert \ \textrm{there\ is} \ (x',y) \in U[k] \ \textrm{with} \ d(x,x') \le K  \}.
\]
For a product set $U=S \times T$, it is more convenient to use the notation $(S,T)[K,k]$ in place of $(S \times T)[K,k]$.
We will refer to the pair $(K,k)$ in the notation $U[K,k]$ as the \textit{enlargement data}.

Let $x_0$ be a chosen fixed point in $X$.
Given a monotone function $h \colon [0, + \infty ) \to [0, + \infty )$, there is a function $h_{x_0} \colon X \to [0, + \infty )$ defined by
\[
h_{x_0} (x) = h (d_X (x_0,x)).
\]

Given two $(X,Y)$-filtered modules $F$ and $G$, an $R$-homomorphism $f \colon F(X \times Y) \to G(X \times Y)$ is \textit{boundedly controlled} if
there are a number $b \ge 0$ and a monotone function $\theta \colon [0, + \infty ) \to [0, + \infty )$ such that
\begin{equation}
fF(U) \subset G(U[b,\theta_{x_0}]) \tag{$\dagger$}
\end{equation}
for all subsets $U \subset X \times Y$ and some choice of $x_0 \in X$.
It is easy to see that this condition is independent of the choice of $x_0$.
If a homomorphism $f$ is boundedly controlled with respect to some choice of parameters 
$b$ and $\theta$, we will say that $f$ is $(b,\theta)$-\textit{controlled}.

The \textit{unrestricted fibred bounded category} $\U_X (Y)$ has $(X,Y)$-filtered modules as objects and the boundedly controlled homomorphisms as morphisms.
Theorem 3.1.6 of \cite{gCbG:18b} shows that $\U_X (Y)$ is a cocomplete semi-abelian category.	

When $Y$ is the one point space, this construction recovers the controlled category $\U (X,R)$ of $X$-filtered $R$-modules used to construct bounded $G$-theory in \cite{gCbG:00} and chapter 2 of \cite{gCbG:18b}.   
In this case, boundedly controlled homomorphisms are characterized by a single parameter $b$, so one can specify that by abbreviating the term to simply $b$-\textit{controlled}.  The construction of an $X$-filtration $F_X$ from a given $(X,Y)$-filtration in Definition \refD{GEPCW2} allows us to view a $(b,\theta)$-controlled homomorphism in $\U_X(Y)$ as a $b$-controlled homomorphism in $\U(X,R)$ via the forgetful functor $T \colon \U_X(Y) \to \U(X,R)$. 	

We now want to restrict to a subcategory of $\U_X (Y)$ that is full on objects with particular properties.  This process consists of two steps that result in a theory with better localization properties.

\begin{DefRef}{StrFib}
An $(X,Y)$-filtered module $F$ is called
\begin{itemize}
\item \textit{split} or $(D,\Delta)$-\textit{split} if there is a number $D \ge 0$ and a monotone function
$\Delta \colon [0, + \infty ) \to [0, + \infty )$ so that
    \[
    F(U_1 \cup U_2) \subset F(U_1 [D,\Delta'_{x_0}]) + F(U_2 [D,\Delta_{x_0}])
    \]
    for each pair of subsets $U_1$ and $U_2$ of $X \times Y$,
\item \textit{lean/split} or $(D,\Delta')$-\textit{lean/split} if there is a number $D \ge 0$ and a monotone function
$\Delta' \colon [0, + \infty ) \to [0, + \infty )$ so that \medskip
\begin{itemize}
\item the $X$-filtered module $F_X$ is $D$-\textit{lean}, in the sense that \[ F_X(S) \subset \sum_{x \in S} F_X(x[D]) \] for every subset $S$ of $X$, while
\item the $(X,Y)$-filtered module $F$ is $(D,\Delta')$-split,
\end{itemize} \medskip
\item \textit{insular} or $(d,\delta)$-\textit{insular} if there is a number $d \ge 0$ and a monotone function
$\delta \colon [0, + \infty ) \to [0, + \infty )$ so that
    \[
    F(U_1) \cap F(U_2) \subset F \big( U_1[d, \delta_{x_0}] \cap U_2[d, \delta_{x_0}] \big)
    \]
    for each pair of subsets $U_1$ and $U_2$ of $X \times Y$.
\end{itemize}

There are two subcategories nested in $\U_X (Y)$.
The category $\mathbf{LS}_X (Y)$ is the full subcategory of $\U_X (Y)$ on objects $F$ that are lean/split and insular.
The category $\B_X (Y)$ is the full subcategory of $\mathbf{LS}_X (Y)$ on objects $F$ such that $F(U)$ is a finitely generated submodule whenever $U \subset X \times Y$ is bounded.  
\end{DefRef}

We proceed to define appropriate exact structures in these categories.
The admissible monomorphisms are precisely the morphisms isomorphic in $\U_X (Y)$ to the filtration-wise monomorphisms and the admissible epimorphisms are those morphisms isomorphic to the filtration-wise epimorphisms.
In other words, the exact structure $\mathcal{E}$ in $\U_X
(Y)$ consists of sequences isomorphic to those
\[
E^{\subdot} \colon \quad E' \xrightarrow{\ i \ } E \xrightarrow{\ j \ } E''
\]
which possess filtration-wise restrictions
\[
E^{\subdot} (U) \colon \quad E' (U) \xrightarrow{\ i \ } E (U) \xrightarrow{\ j \ } E'' (U)
\]
for all subsets $U \subset (X,Y)$, and each $E^{\subdot} (U)$ is an exact sequence of $R$-modules.

Both $\mathbf{LS}_X (Y)$ and $\B_X (Y)$ 
are closed under extensions in $\U_X (Y)$. 
Therefore, they are themselves exact categories,
and the inclusion
$\mathcal{B}_X (Y) \to \B_X (Y)$
is an exact embedding, as we projected.

\medskip

There is a useful invariant of a finitely generated group $\Gamma$ that is defined in terms of the exact category $\B_{\Gamma} (\point)$ in \cite{gCbG:15}.  Here $\Gamma$ can be given the word metric associated to any of the finite generating sets.  The left multiplication action gives an action of $\Gamma$ on $\B_{\Gamma} (\point)$.

Recall that Theorem \refT{JGun} provides an interpretation to the $K$-theory of a group ring $R[\Gamma]$ in terms of relative homotopy fixed points of the additive category $\mathcal{B} (X,R)$, which can be viewed as $\mathcal{B}_{\Gamma} (\point)$.

\begin{ExRefName}{HFPGG}{Bounded $G$-theory of a finitely generated group} 
In the case where $\mathcal{C}$ is the exact category $\B_{\Gamma} (\point)$ and $\mathcal{M}$ is the subcategory of the filtration zero morphisms, the \textit{bounded $G$-theory of $\Gamma$} is defined to be the nonconnective $K$-theory of the relative homotopy fixed points $\B_{\Gamma} (\point)^{h\Gamma}$, denoted $\Gnc (R[\Gamma])$.    
\end{ExRefName}

Notice that this definition makes sense even when the group ring is not Noetherian unlike the much more restrictive situation with the usual $G$-theory defined only for Noetherian rings.  

\begin{ThmRef}{JGun2}
There is an exact subcategory of finitely generated $\Gamma$-modules for an arbitrary finitely generated group $\Gamma$ such that its relative homotopy fixed points have Quillen $K$-theory with features similar to $G$-theory of group rings. In particular, it has a Cartan map from the $K$-theory of $R[\Gamma]$.
\end{ThmRef}

\begin{proof}
The category is equivalent to $\B_{\Gamma} (\point)$.  We refer to sections 2 and 3 of \cite{gCbG:15} for details.
The clear resemblance to Definition \refD{NMEWQS} and the identification of $\mathcal{B}^{\Gamma, 0} (X,R)^{\Gamma}$ with $\mathcal{B}_{\Gamma} (\point)^{h\Gamma}$ allow us to induce the Cartan map $\Knc (R[\Gamma]) \to \Gnc (R[\Gamma])$ from the exact inclusion $\mathcal{B}_{\Gamma} (\point) \to \B_{\Gamma} (\point)$ above.
\end{proof}

Suppose $C$ is a subset of $Y$.  Let $\B_X (Y)_{<C}$ be the full subcategory of $\B_X (Y)$ on objects $F$ such that
    \[
    F (X, Y) \subset F \big( (X, C)[r,\rho_{x_0}]) \big)
    \]
    for some number $r \ge 0$ and an order preserving function $\rho \colon [0,+\infty) \to [0,+\infty)$.
    
Recall that a \textit{Serre subcategory} of
an exact category is a full subcategory which is closed under
exact extensions and closed under passage to admissible subobjects and
admissible quotients.
Proposition 3.3.3 of \cite{gCbG:18b} verifies that
$\B_X (Y)_{<C}$ is a Serre subcategory of $\B_X (Y)$.

The second step in restricting to subcategories with good localization properties is done via introducing a structure called a grading.

Given an arbitrary $R$-submodule $F'$ of $F$ in $\U_X(Y)$, we can assign to $F'$ the \textit{standard $(X,Y)$-filtration} $F' (U)=F(U)\cap F'$.

Let $\mathcal{M}^{\ge 0}$ be the set of all monotone functions $\delta \colon [0, + \infty) \to [0, + \infty)$.

Let $\mathcal{P}_X (Y)$ be the subcategory of $\mathcal{P}(X,Y)$ consisting of
all subsets of the form $(X,C)[D,\delta_{x_0}]$ for some choices of a subset $C \subset Y$, a number $D \ge 0$, and a function $\delta \in \mathcal{M}^{\ge 0}$.

\begin{DefRef}{LScov22}
Given an object $F$ of $\B_X (Y)$, a $Y$-\textit{grading} of $F$
is a functor
$\mathcal{F} \colon \mathcal{P}_X (Y) \to  \mathcal{I}(F)$
with the following properties:
\begin{itemize}
\item the submodule $\mathcal{F} ((X,C)[D,\delta_{x_0}])$, with the standard $(X,Y)$-filtration induced from $F$, is an object of $\B_X (Y)$,
\item there is an enlargement data $(K,k)$ such that
\[
F((X,C)[D,\delta_{x_0}]) \subset \mathcal{F} ((X,C)[D,\delta_{x_0}]) \subset F((X,C) [D + K, \delta_{x_0} + k_{x_0}]),
\]
for all subsets in $\mathcal{P}_X (Y)$.
\end{itemize}
\end{DefRef}

We say that an object $F$ of $\B_X (Y)$ is $Y$-\textit{graded} if there exists a $Y$-grading of $F$, but the grading itself is not specified, and define $\G_X (Y)$ as the full subcategory of $\B_X (Y)$ on $Y$-graded filtered modules.  

We will summarize some additional required results from section 3.4 of \cite{gCbG:18b}.

\begin{ThmRef}{VFGHJ2}
The subcategory $\G_X (Y)$ is closed under both isomorphisms and exact extensions in $\B_X (Y)$.
Therefore, $\G_X (Y)$ is an exact subcategory of $\B_X (Y)$.

The restriction to $Y$-gradings in $\B_X (Y)_{<C}$ gives a full exact subcategory $\G_X (Y)_{<C}$ which is a Serre subcategory of $\G_X (Y)$.

Given a graded object $F$ in $\G_X (Y)$,
we assume that $F$ is $(D,\Delta')$-split and $(d, \delta)$-insular and is graded by $\mathcal{F}$.  For a subset $U$ from the family $\mathcal{P}_X (Y)$,
the submodule $\mathcal{F} (U)$ has the following properties:
\begin{enumerate}
\item $\mathcal{F} (U)$ is graded by $\mathcal{F}_U (T) = \mathcal{F} (U) \cap \mathcal{F} (T)$,
\item $F(U) \subset \mathcal{F} (U) \subset F(U[K,k])$ for some fixed enlargement data $(K,k)$,
\item if $q \colon F \to H$ is the quotient of the inclusion $i \colon \mathcal{F} (U) \to F$ and $F$ is $(D, \Delta')$-lean/split, then $H$ is supported on $(X \setminus U) [2D,2\Delta']$,
\item $H (U[-2D -2d, -2\Delta' - 2\delta]) =0$.
\end{enumerate}
\end{ThmRef}

We assume that the reader is familiar with Quillen $K$-theory of exact categories.  This theory can be applied to both $\G_X (Y)$ and $\G_X (Y)_{<C}$.  The result can be viewed as the spectra $G_X (Y)$ and $G_X (Y)_{<C}$.  The stable homotopy groups of these spectra are the Quillen $K$-groups of the exact categories.

Finally, the main goal of this section is a homotopy fibration
\[
G_X (Y)_{<C} \longrightarrow G_X (Y) \longrightarrow G_X (Y,C)
\]
where $G_X (Y,C)$ is the $K$-theory of a certain quotient category $\G_X (Y)/\G_X (Y)_{<C}$. 

For simplicity we will use the notation $\G$ for $\G_X (Y)$ and $\C$ for the Serre subcategory $\G_X (Y)_{<C}$ of $\G$.

There is a class
of \textit{weak equivalences} $\Sigma (C)$ in $\G$ which consist of all finite
compositions of admissible monomorphisms with cokernels in $\C$ and
admissible epimorphisms with kernels in $\C$.  We need the class $\Sigma (C)$ to admit calculus of right fractions.  This follows from  \cite[Lemma 1.13]{mS:03} and the fact that $\C$ in $\G$ is right
filtering, in the sense that each morphism $f \colon F_1 \to F_2$ in $\G$, where $F_2$
is an object of $\C$, factors through an admissible epimorphism $e
\colon F_1 \to \overline{F}_2$, where $\overline{F}_2$ is in $\C$.  The latter fact is \cite[Lemma 3.5.6]{gCbG:18b}.

The category $\mathbf{G}/\mathbf{C}$ is the localization $\G [\Sigma (C)^{-1}]$.
From \cite[Theorem 3.5.8]{gCbG:18b}, it is an
exact category where the short sequences are
isomorphic to images of exact sequences from $\G$.

There is an intrinsic reformulation of the homotopy fibration because the essential full image of the evident inclusion of $\G_X (C)$ in $\G$ is precisely $\C$.  This gives a homotopy fibration
\begin{gather*}
{G}_X (C) \longrightarrow {G}_X (Y) \longrightarrow {G}_X (Y,C).
\end{gather*}

One quick consequence is the ability to relativize the old constructions.  If $Y'$ is any subset of $Y$, one obtains the relative theory $G_X (Y,Y')$.

Another easy application is a nonconnective delooping that applies to  all of the theories we have defined.  For example in the basic case,   
\[
\GncX (Y)  \, = \,  \hocolim{k>0}
\Omega^{k} {G}_{X} (Y \times \mathbb{R}^k).
\]
This uses the standard Eilenberg swindle trick and can be seen in section 4.2 of \cite{gCbG:18b}.

\SecRef{Fibrewise excision in equivariant fibred $G$-theory}{TTITS}

It is well-known that Quillen $K$-theory of an exact category can be obtained equivalently as Waldhausen's $K$-theory of bounded chain complexes in the category.  The cofibrations are then the chain maps which are the degree-wise admissible monomorphisms. The weak equivalences are the chain maps whose mapping cones are homotopy equivalent to acyclic complexes. An exposition with a number of details verified specifically for bounded $G$-theory can be found in \cite[section 4]{gCbG:00}.  The Waldhausen theory setting is crucial in proving the excision theorem in that the Approximation Theorem \cite[Theorem 4.5]{gCbG:00} becomes essential.  We will indicate passage from an exact category to the derived category of bounded chain complexes by prefixing ``$\mathrm{ch}$'' in front of the name of the exact category.

We proceed to define the equivariant fibred $G$-theory.

The basic setting consists of 
\begin{itemize}
	\item two proper metric spaces $X$ and $Y$,
	\item an arbitrary subset $Y'$ of $Y$,
	\item a $\Gamma$-action on $X$ by isometries, and 
	\item a bounded action of $\Gamma$ on $Y$.  This is an action such that for each $\gamma$ in $\Gamma$ the set of real numbers $W_{\gamma} = \{ d(x, \gamma(x)) \}$ is bounded from above.
\end{itemize}

 \begin{RemRef}{EQF} In a number of situations, we will be specifying subcategories closed under the $\Gamma$-action by subsets that are arbitrary, therefore certainly not closed under the action.  This works due to the boundedness of the action.  For example, if we have a subset $C \subseteq X$ and define a subcategory as the set of modules supported on some neighborhood of $C$, then this subcategory is closed under the $\Gamma$-action provided the action is bounded.  This would definitely not hold were the action not bounded. \end{RemRef} 

Consider the exact category $\G_{\Gamma} (Y)$ with the induced action by $\Gamma$, in the case $X$ is the group $\Gamma$ with a word metric, acting on itself by isometries via the left multiplication. Since the action on $Y$ is bounded, we have the quotient exact category $\G_{\Gamma} (Y, Y')$.  

\begin{NotRef}{MIYB}
	If $Z$ is another arbitrary subset of $Y$, it is also useful to consider the full exact subcategory $\G_{\Gamma} (Y, Y')_{<Z}$ which we denote $\G_{\Gamma} (Y, Y',Z)$.
\end{NotRef}

Recasting the definitions from \refD{RLL2}, we define the Waldhausen category $\mathcal{G}_{\Gamma,0} (Y, Y',Z)$ to be the full subcategory of $\Fun(\EGamma,\mathrm{ch}\G_{\Gamma} (Y, Y',Z))$ on those functors that send morphisms in $\EGamma$ to degree-wise $0$-controlled homomorphisms of $\Gamma$-filtered modules.  

\begin{DefRef}{HDJBF}
The \textit{equivariant fibred $G$-theory} is
\[
G^{\Gamma} (Y,Y',Z) = \Omega K (\vert wS. \mathcal{G}_{\Gamma,0} (Y, Y',Z) \vert). 
\]
This is a functor from the category of triples $(Y,Y',Z)$, where both $Y'$ and $Z$ are subspaces of $Y$ but not necessarily subspaces of each other, and uniformly expansive maps of triples to the category of spectra.
\end{DefRef}

Now we turn to the construction of fibred homotopy fixed points.
There is a forget control functor $T \colon \G_X (Y,Y',Z) \to \U_X (Y,Y',Z)$  sending $F$ to $F_X$.  
Since $\Gamma$ acts on $X$ by isometries, it also acts on $\U_X (Y,Y',Z)$.   The combination of this action and a bounded action on $Y$ induces an action on $\G_X (Y,Y',Z)$.  With these choices, $T$ is an equivariant functor.

\begin{DefRefName}{RLL5}{Fibred homotopy fixed points in bounded $G$-theory}
This is a special case of a relative homotopy fixed points, as defined in \refD{RLL2}, with the choices of $\mathcal{C}$ and $\mathcal{M}$ as follows.
\begin{itemize}
	\item[--] the category $\mathcal{C}$ is $\G_X (Y,Y',Z)$, 
    \item[--] the subcategory $\mathcal{M}$ consists of all controlled morphisms $\phi$ in $\mathcal{C}$ with the property that $T(\phi)$ is bounded by $0$ as homomorphisms controlled over $X$.
\end{itemize}

Let us recapitulate what this definition entails in the case $X$ is the group $\Gamma$ with a word metric.

The \textit{fibred homotopy fixed points of a triple $(Y,Y^{\prime}, Z)$} is the category $\G^{h\Gamma} (Y, Y',Z)$ with objects which are sets of data $( \{ F_{\gamma} \},\{ \psi_{\gamma} \} )$ where
\begin{itemize}
\item $F_{\gamma}$ is an object of $\G_{\Gamma} (Y,Y',Z)$ for each $\gamma$ in $\Gamma$,
\item $\psi_{\gamma}$ is an isomorphism $F_e \to F_{\gamma}$ in $\G_{\Gamma} (Y,Y',Z)$,
\item $\psi_{\gamma}$ is $0$-controlled when viewed as a morphism in $\U_{\Gamma} (Y,Y',Z)$,
\item $\psi_e = \id$,
\item $\psi_{\gamma_1 \gamma_2} = \gamma_1 \psi_{\gamma_2} \circ \psi_{\gamma_1}$ for all $\gamma_1$, $\gamma_2$ in $\Gamma$.
\end{itemize}
The morphisms $( \{ F_{\gamma} \},\{ \psi_{\gamma} \} ) \to ( \{ F'_{\gamma} \},\{ \psi'_{\gamma} \} )$
are collections of morphisms $\phi_{\gamma} \colon F_{\gamma} \to F'_{\gamma}$ in $\G_X (Y,Y',Z)$ such that the squares
\[
\xymatrix{
F_e \ar[rr]^-{\psi_{\gamma}} \ar[dd]_-{\phi_e} &&F_{\gamma} \ar[dd]^-{\phi_{\gamma}}\\
\\
F'_e \ar[rr]^-{\psi'_{\gamma}} &&F'_{\gamma}
}
\]
commute for all $\gamma$.

The exact structure on $\G^{h\Gamma} (Y, Y',Z)$ is induced from that on $\G_{\Gamma} (Y,Y',Z)$ as follows.
A morphism $\phi$ in $\G^{h\Gamma} (Y, Y',Z)$ is an admissible monomorphism if $\phi_e \colon F \to F'$ is an admissible monomorphism in $\G_{\Gamma} (Y, Y',Z)$.  This of course implies that all structure maps $\phi_{\gamma}$ are admissible monomorphisms.  Similarly, a morphism $\phi$ is an admissible epimorphism if $\phi_e \colon F \to F'$ is an admissible epimorphism.  This gives $\G^{h\Gamma} (Y, Y',Z)$ an exact structure.
\end{DefRefName}

Since the induced $\Gamma$-action on $S. \mathcal{G}_{\Gamma,0} (Y, Y',Z)$ commutes with taking fixed points, we have the following fact.

\begin{PropRef}{HXCV}
	The fixed point spectrum $G^{\Gamma} (Y,Y',Z)^{\Gamma}$ is equivalent to the $K$-theory of the relative homotopy fixed point category $\G^{h\Gamma} (Y,Y',Z)$.
\end{PropRef}

We proceed to consider multiple bounded actions of $\Gamma$ on $Y$.  Let $\beta (Y)$ be the set of all such actions.  
Let $\mathcal{F}$ be the functor that assigns to a set $Z$ the partially ordered set of finite subsets of $Z$.

\begin{DefRef}{MPEW}
	For any $S$ in $\mathcal{F}(\beta (Y))$ we define $Y_S$ as the metric space which is the disjoint union $\bigsqcup_{s \in S} Y_s$, where $Y_s$ are copies of $Y $ with the specified action.  The metric on $Y_S$ is induced by the requirement that  it restricts to the metric from $Y$ in each $Y_s$ and for the same point $y$ in different components the distance $d (y_s, y_{s'})$ equals 1.
\end{DefRef}
    
    Clearly, the action of $\Gamma$ on $Y_S$ is bounded.  

As a consequence of Proposition \refP{HXCV}, for each choice of finite subset $S$ of $\beta (Y)$, the spectrum $G^{\Gamma} (Y_S, Y'_S,Z)^{\Gamma}$ is the Quillen $K$-theory spectrum of $\G^{h\Gamma} (Y_S, Y'_S,Z)$, where $Z$ is a subset of $Y_S$. 

\begin{ThmRef}{FibrG}
Let $C$ be an arbitrary subset of $Y$.
There is a homotopy fibration
\begin{gather*}
G^{\Gamma} (Y_S, Y'_S,Z)^{\Gamma}_{<C} \longrightarrow G^{\Gamma} (Y_S, Y'_S,Z)^{\Gamma} \longrightarrow G^{\Gamma} (Y_S, Y'_S,Z)^{\Gamma}_{>C},
\end{gather*}
where $G^{\Gamma} (Y_S, Y'_S,Z)^{\Gamma}_{>C}$ stands for $K$-theory of the exact quotient $\G^{h\Gamma} (Y_S, Y'_S,Z)_{>C}$ of $\G^{h\Gamma} (Y_S, Y'_S,Z)$ by $\G^{h\Gamma} (Y_S, Y'_S,Z)_{<C}$.
In the absolute case, there is an equivalence $\G^{h\Gamma}(Y)_{<C} \simeq \G^{h\Gamma} (C)$, and so there a special case of the homotopy fibration
\begin{gather*}
G^{\Gamma} (C)^{\Gamma} \longrightarrow G^{\Gamma} (Y)^{\Gamma} \longrightarrow G^{\Gamma} (Y)^{\Gamma}_{>C}.
\end{gather*}
\end{ThmRef}

\begin{proof}  
	Taking into account Remark \refR{EQF}, the fact that $\G_{\Gamma} (Y_S, Y'_S,Z)_{<C}$ is an idempotent complete Serre subcategory of $\G_{\Gamma} (Y_S, Y'_S,Z)$ implies immediately that $\G^{h\Gamma} (Y_S, Y'_S,Z)^{h\Gamma}_{<C}$ is an idempotent complete Serre subcategory of $\G^{h\Gamma} (Y_S, Y'_S,Z)$.  
	The main technical result of Schlichting \cite{mS:03} is a fibration theorem which requires $\G^{h\Gamma} (Y_S, Y'_S,Z)_{<C}$ to satisfy two additional properties, right filtering and right $s$-filtering. Both of these properties follow directly from the estimates in the proofs of Lemma 3.5.6 and Theorem 3.5.8 in \cite{gCbG:18b}. 
\end{proof}

Our first application of the fibration is to deloop $G^{\Gamma} (Y_S, Y'_S,Z)^{\Gamma}$ and related spectra following the strategy after Pedersen/Weibel \cite{ePcW:85}.  

Let $\mathbb{R}$, $\mathbb{R}^{\ge 0}$, and $\mathbb{R}^{\le 0}$
denote the metric spaces of the reals, the nonnegative reals, and
the nonpositive reals with the Euclidean metric.  Then there is the following
map of homotopy fibrations
\[
\xymatrix{
G^{\Gamma} (Y_S)^{\Gamma} \ar[r] \ar[d] &G^{\Gamma} (Y_S \times \mathbb{R}^{\ge 0})^{\Gamma} \ar[r] \ar[d] &G^{\Gamma} (Y_S \times \mathbb{R}^{\ge 0})^{\Gamma}_{>Y_S \times 0} \ar[d]^-{K(I)}  \\
G^{\Gamma} (Y_S \times \mathbb{R}^{\le 0})^{\Gamma} \ar[r] &G^{\Gamma} (Y_S \times \mathbb{R})^{\Gamma} \ar[r] &G^{\Gamma} (Y_S \times \mathbb{R})^{\Gamma}_{> Y_S \times \mathbb{R}^{\le 0}}
}
\]
The map $K(I)$ is induced by the inclusion $I$ of the quotient categories. 

\begin{ThmRef}{HQASW}
	$K(I)$ is a weak equivalence of connective spectra. 
\end{ThmRef}

\begin{proof}
This follows from the Approximation Theorem applied to $I$. 
The first condition of the theorem is evident. 
To check the second condition, consider a chain complex $F^{\subdot}$ in 
$\G^{h\Gamma} (Y_S \times \mathbb{R}^{\ge 0})_{>Y_S \times 0}$.
By the
nature of the objects and the explanation in Remark \refR{EQF}, all maps in $F^{\subdot}$ and their control features are determined by the values on the objects $F^i_e$ of $\G_{\Gamma} (Y_S \times \mathbb{R})$.
So we can specify $F^{\subdot}$ by the chain complex
\[
F^{\subdot}_e \colon \quad 0 \longrightarrow F^1_e \xrightarrow{\
\phi_1\ } F^2_e \xrightarrow{\ \phi_2\ }\ \dots\
\xrightarrow{\ \phi_{n-1}\ } F^n_e \longrightarrow 0
\]
in $\G_{\Gamma} (Y_S \times \mathbb{R}^{\ge 0})_{>Y_S \times 0}$.
Given a chain complex $G^{\subdot}$ in $\G^{h\Gamma} (Y_S \times \mathbb{R})_{>Y_S \times \mathbb{R}^{\le 0}}$, 
we can apply the same reasoning to $G^{\subdot}$. 
Now a chain map $g \colon F^{\subdot} \to G^{\subdot}$ is given uniquely by a chain map $F^{\subdot}_e \to G^{\subdot}_e$ where $G^{\subdot}_e$ is the chain complex
\[
0 \longrightarrow G^1_e \xrightarrow{\
\psi_1\ } G^2_e \xrightarrow{\ \psi_2\ }\ \dots\
\xrightarrow{\ \psi_{n-1}\ } G^n_e \longrightarrow 0
\]
in $\G_{\Gamma} (Y_S \times \mathbb{R})_{>Y_S \times \mathbb{R}^{\le 0}}$.
Also observe that if $F^i_e$ is supported on a neighborhood of $C \subset Y$ then so are all of $F^i_{\gamma}$.  
This allows us to transport from \cite{gCbG:18b} the rest of the argument for a non-equivariant Lemma 4.2.4.
Alternatively, we can refer to the end of the proof of Theorem \refT{Exc2} where the details are spelled out in even greater generality.
\end{proof}

The spectra $G^{\Gamma} (Y_S \times \mathbb{R}^{\ge 0})^{\Gamma}$ and $G^{\Gamma} (Y_S \times \mathbb{R}^{\le 0})^{\Gamma}$ are contractible as $K$-theory spectra of flasque categories.  This is the standard consequence of the shift functor $T$ in the positive
(respectively negative)
direction along $\mathbb{R}^{\ge 0}$ (respectively $\mathbb{R}^{\le 0}$) interpreted as an exact endofunctor.  A natural equivalence of functors 
$1 \oplus T \cong T$ and the
Additivity Theorem give contractibility.

From the map of fibrations, we obtain a map of spectra $G^{\Gamma} (Y_S)^{\Gamma} \to \Omega G^{\Gamma} (Y_S \times \mathbb{R})^{\Gamma}$
which induces isomorphisms of $K$-groups
in positive dimensions.
Iterating this construction for $k \ge 2$ gives weak equivalences
\[
\Omega^k G^{\Gamma} (Y_S \times \mathbb{R}^k)^{\Gamma}
\longrightarrow \Omega ^{k+1}
G^{\Gamma} (Y_S \times \mathbb{R}^{k+1})^{\Gamma}.
\]

\begin{DefRef}{RealExQ}
The nonconnective delooping of algebraic $K$-theory of the fibred homotopy fixed points is the spectrum
\[
\widetilde{G}^{\Gamma} (Y_S)^{\Gamma}  \, = \,  \hocolim{k>0}
\Omega^{k} G^{\Gamma} (Y_S \times \mathbb{R}^k)^{\Gamma}.
\]
In the case $Y$ is the one point space, $\widetilde{G}^{\Gamma} (Y)^{\Gamma}$ coincides with the nonconnective $G$-theory of the group ring $R[\Gamma]$ defined by the authors in \cite{gCbG:15}.
\end{DefRef}

The discussion leading up to Definition \refD{RealExQ} can be repeated verbatim for other Serre subcategory pairs.  For example, the subcategory $\G^{\Gamma} (Y_S \times \mathbb{R}^k)^{\Gamma}_{<C \times \mathbb{R}^k}$ is evidently a Serre subcategory of $\G^{\Gamma} (Y_S \times \mathbb{R}^k)^{\Gamma}$ for any choice of subset $C \subset Y_S$.
We define
\[
\widetilde{G}^{\Gamma} (Y_S)^{\Gamma}_{<C}  \, = \,  \hocolim{k>0}
\Omega^{k} {G}^{\Gamma} (Y_S \times \mathbb{R}^k)^{\Gamma}_{<C \times \mathbb{R}^k}
\]
and
\[
\widetilde{G}^{\Gamma} (Y_S)^{\Gamma}_{<C_1, C_2}  \, = \,  \hocolim{k>0}
\Omega^{k} {G}^{\Gamma} (Y_S \times \mathbb{R}^k)^{\Gamma}_{<C_1 \times \mathbb{R}^k, \, C_2 \times \mathbb{R}^k}.
\]

\begin{DefRef}{GEPCWpr}
Let $Y'$, $Y_1$ and $Y_2$ be arbitrary subsets of $Y$ so that $Y_1$ and $Y_2$ form a covering of $Y$.
There are corresponding subsets $Y'_S$, $Y_{1,S}$ and $Y_{2,S}$ of $Y_S$ obtained as $Y'_S = \bigsqcup Y'_s$, $Y_{1,S} = \bigsqcup Y'_{1,s}$ and $Y_{2,S} = \bigsqcup Y'_{2,s}$.
It is now straightforward to define nonconnective spectra
\[
\widetilde{G}^{\Gamma} (Y_S,Y'_S)^{\Gamma}  \, = \,  \hocolim{k>0}
\Omega^{k} {G}^{\Gamma} (Y_S \times \mathbb{R}^k, Y'_S \times \mathbb{R}^k)^{\Gamma},
\]
\[
\widetilde{G}^{\Gamma} (Y_S,Y'_S)^{\Gamma}_{<Y_i}  \, = \,  \hocolim{k>0}
\Omega^{k} {G}^{\Gamma} (Y_S \times \mathbb{R}^k, Y'_S \times \mathbb{R}^k)^{\Gamma}_{<Y_{i,S} \times \mathbb{R}^k},
\]
\[
\widetilde{G}^{\Gamma} (Y_S,Y'_S)^{\Gamma}_{<Y_1, Y_2}  \, = \,  \hocolim{k>0}
\Omega^{k} {G}^{\Gamma}  (Y_S \times \mathbb{R}^k, Y'_S \times \mathbb{R}^k)^{\Gamma}_{<Y_{1,S} \times \mathbb{R}^k, \, Y_{2,S} \times \mathbb{R}^k}.
\]
\end{DefRef}

\begin{ThmRef}{Exc2}
Suppose $Y_1$ and $Y_2$ are subsets of a metric space $Y$, and $Y = Y_1 \cup Y_2$.
There is a homotopy pushout diagram of spectra
\[
\xymatrix{
 \widetilde{G}^{\Gamma}(Y_S,Y'_S)^{\Gamma}_{<Y_1,Y_2} \ar[r] \ar[d]
&\widetilde{G}^{\Gamma}(Y_S,Y'_S)^{\Gamma}_{<Y_1} \ar[d] \\
 \widetilde{G}^{\Gamma} (Y_S,Y'_S)^{\Gamma}_{<Y_2} \ar[r]
&\widetilde{G}^{\Gamma}(Y_S,Y'_S)^{\Gamma}
}
\]
where the maps of spectra are induced from the exact inclusions.
If we define
\[
{E}^{\Gamma} (Y,Y') = \hocolim{U \in \mathcal{F}(\beta (Y))} \widetilde{G}^{\Gamma} (Y_S, Y'_S)^{\Gamma},
\]
the excision theorem also holds on the level of ${E}^{\Gamma}$.  There is a homotopy pushout diagram of spectra
\[
\xymatrix{
 {E}^{\Gamma} (Y,Y')_{<Y_1,Y_2} \ar[r] \ar[d]
&{E}^{\Gamma} (Y,Y')_{<Y_1} \ar[d] \\
 {E}^{\Gamma} (Y,Y')_{<Y_2} \ar[r]
&{E}^{\Gamma}(Y,Y')
}
\]
\end{ThmRef}

\begin{proof}
There is a homotopy pushout
\[
\xymatrix{
 \widetilde{G}^{\Gamma} (Y_S,Y'_S)^{\Gamma}_{< Y_1, Y_2} \ar[r] \ar[d]
&\widetilde{G}^{\Gamma} (Y_S,Y'_S)^{\Gamma}_{<Y_1} \ar[d] \\
 \widetilde{G}^{\Gamma} (Y_S,Y'_S)^{\Gamma}_{<Y_2} \ar[r]
&\widetilde{G}^{\Gamma} (Y_S,Y'_S)^{\Gamma}
}
\]
obtained from the map of the fibration sequences
\[
\xymatrix{
\widetilde{G}^{\Gamma} (Y_S,Y'_S)^{\Gamma}_{< Y_1, Y_2} \ar[r] \ar[d] 
&\widetilde{G}^{\Gamma} (Y_S,Y'_S)^{\Gamma}_{<Y_1} \ar[r] \ar[d] 
&\widetilde{G}^{\Gamma} (Y_S,Y'_S)^{\Gamma}_{<Y_1, >Y_2} \ar[d] \\ 
\widetilde{G}^{\Gamma} (Y_S,Y'_S)^{\Gamma}_{<Y_2} \ar[r] 
&\widetilde{G}^{\Gamma} (Y_S,Y'_S)^{\Gamma} \ar[r] 
&\widetilde{G}^{\Gamma} (Y_S,Y'_S)^{\Gamma}_{>Y_2}
}
\]
both obtained from Theorem \refT{FibrG}.
The map 
\[
\widetilde{G}^{\Gamma} (Y_S,Y'_S)^{\Gamma}_{<Y_1, >Y_2} \longrightarrow \widetilde{G}^{\Gamma} (Y_S,Y'_S)^{\Gamma}_{>Y_2}
\]
induced from the exact inclusion $J \colon \G^{h\Gamma} (Y_S,Y'_S)_{<Y_1, >Y_2} \to \G^{h\Gamma}  (Y_S,Y'_S)_{>Y_2}$
is again an equivalence.  It should be instructive to spell out the crucial application of the Approximation Theorem.
Consider a chain complex $F^{\subdot}$ in 
$\G^{h\Gamma} (Y_S,Y'_S)_{<Y_1, >Y_2}$.
All maps in $F^{\subdot}$ and their control features are determined by the values on the objects 
$F^i_e$ of $\G_{\Gamma} (Y_S,Y'_S)$, so $F^{\subdot}$ can be given by the chain complex
\[
F^{\subdot}_e \colon \quad 0 \longrightarrow F^1_e \xrightarrow{\
\phi_1\ } F^2_e \xrightarrow{\ \phi_2\ }\ \dots\
\xrightarrow{\ \phi_{n-1}\ } F^n_e \longrightarrow 0
\]
in $\G_{\Gamma} (Y_S,Y'_S)_{<Y_1, >Y_2}$.
Applying the same reasoning to a chain complex $G^{\subdot}$ in $\G^{h\Gamma} (Y_S,Y'_S)_{>Y_2}$, let $G^{\subdot}_e$ be the chain complex
\[
0 \longrightarrow G^1_e \xrightarrow{\
\psi_1\ } G^2_e \xrightarrow{\ \psi_2\ }\ \dots\
\xrightarrow{\ \psi_{n-1}\ } G^n_e \longrightarrow 0
\]
in $\G_{\Gamma} (Y_S,Y'_S)_{>Y_2}$.
A chain map $g \colon F^{\subdot} \to G^{\subdot}$ can be given by a chain map $g' \colon F^{\subdot}_e \to G^{\subdot}_e$. 
Since the action is bounded, if $F^i_e$ is supported near a neighborhood of $Y_{1,S} \subset Y_S$ then so are all $F^i_{\gamma}$. 
 
Suppose all $F^i_e$ and $G^i_e$ are $(D,\Delta')$-lean/split and $(d,\delta)$-insular, and there is a number $r \ge 0$ and a monotone function $\rho \colon [0,+\infty) \to [0,+\infty)$ such that there are containments 
$F^i_e \subset F^i_e ((\Gamma,Y_{1,S}) [r,\rho_{x_0}])$
for all $0 \le i \le n$.
Suppose also that $(b,\theta)$ can be used as bounded control data for all maps $\phi_i$, $\psi_i$, and $g'_i$.
Suppose also that $(K,k)$ is an enlargement data for the chosen grading of $G_e$.
We define a submodule
$F_e^{\prime{i}}$ as the submodule $\mathcal{G}_e^{i} ((\Gamma,Y_{1,S}) [r+3ib,\rho_{x_0} + 3i\theta_{x_0}])$ in the chosen grading of $G_e$
and define $\xi_{i} \colon F^{\prime{i}}_e \to F^{\prime{i+1}}_e$ to be the restrictions of
$\psi_{i}$ to $F^{\prime{i}}_e$. This gives a chain subcomplex
$(F^{\prime{i}},\xi_{i})$ of $(G^{i},\psi_{i})$ in $\G_{\Gamma} (Y_S,Y'_S)$ with the inclusion $i \colon F^{\prime{i}} \to G^{i}$.  Notice that we have the induced chain map $\overline{g} \colon F^{\subdot} \to F^{\prime\subdot}$ in $\G_{\Gamma} (Y_S,Y'_S)_{<Y_1}$ so that $g = i J(\overline{g})$.
It remains to prove that the cokernel $C^{\subdot}$ of $i$ is in $\G_{\Gamma} (Y_S,Y'_S)_{<Y_2}$.
Since
\[
F^{\prime{i}}_e \subset {G}^{i}_e ((\Gamma,Y_{1,S}) [r+3ib +K,\rho_{x_0} + 3i\theta_{x_0} +k_{x_0}]),
\]
each $C^i$ is supported on
\begin{multline*}
(\Gamma,Y_S \setminus Y_{1,S})[2D +2d - r - 3ib - K, 2\Delta'_{x_0} + 2\delta_{x_0} - \rho_{x_0} - 3i\theta_{x_0} - k_{x_0}] \\
\subset
(\Gamma, Y_{2,S}) [2D +2d, 2\Delta'_{x_0} + 2\delta_{x_0}],
\end{multline*}
This shows that the complex $C^{\subdot}$ is indeed in $\G_{\Gamma} (Y_S,Y'_S)_{<Y_2}$. 
\end{proof}

\textbf{A technical remark.}
This observation is warranted as we contrast Theorem~\refT{Exc2} and its proof with inability to use other, more standard methods in bounded algebra based on Karoubi filtrations in order to prove similar facts in $K$-theory.  The key idea in the proof is still the commutative diagram from Cardenas/Pedersen \cite[section 8]{mCeP:97} transported from bounded $K$-theory to fibred $G$-theory.  Cardenas and Pedersen use Karoubi quotients and the Karoubi fibrations in order to establish their diagram.  One of the crucial points in \cite{mCeP:97} is that the functor $I$ between the Karoubi quotients is an isomorphism of categories.
In fibred $G$-theory the situation is more complicated: $I$ is not necessarily full and, therefore, not an isomorphism of categories.  However we can see here just as in the analogous Theorem 4.4.2 in \cite{gCbG:18b}, the Approximation Theorem suffices to prove that $K(I)$ is nevertheless a weak equivalence.

We make an explicit statement that does not hold in $K$-theory.  Let $V_{\Gamma} (Y)$ be the $K$-theory of the fibred homotopy fixed points $\mathcal{B}^{w\Gamma} (Y)$ defined in \refN{NMDKV}.  In these terms, we don't know whether, or under what conditions on $\Gamma$ and $Y$, 
\[
\hocolim{U \in \mathcal{U}} V^{\Gamma} (U) \longrightarrow V^{\Gamma} (Y)
\]
is an equivalence in the context of Theorem \refT{Exc2}.
Now Example 5.2 in \cite{gCbG:18a} demonstrates in the most basic geometric situation that $\mathcal{B}^{w\Gamma} (Y)$ fails to be Karoubi filtered by the natural choice of subcategory.  Through indirect ways related to the work on the Borel conjecture sketched in section \refS{FFLLIA}, we know that Karoubi filtrations should be impossible to use to compute fibred homotopy fixed points in full generality because of the well-known counterexamples to the isomorphism conjecture for the assembly map in cases of non-regular rings $R$.

\medskip

Suppose $\mathcal{U}$ is a finite covering of $Y$ that is closed under intersections and such that the family of all subsets $U$ in $\mathcal{U}$ together with $Y'$ are pairwise coarsely antithetic.
The extra conditions in the second statement ensure that the covering is in fact by complete representatives of a covering by ``coarse families'' in the language introduced in \cite[section 4.3]{gCbG:18b}.  

We define the homotopy colimit
\[
\mathcal{E}^{\Gamma} (Y,Y')_{< \mathcal{U}} \, = \, \hocolim{U \in \mathcal{U}} E^{\Gamma} (Y, Y')_{<U}.
\]

All of the above discussion can be restricted to full subcategories on objects supported near an arbitrary subset $Z$ of $Y_S$, so we have the following general statement.

\begin{ThmRefName}{PPPOI8fin}{Fibrewise bounded excision}
The natural map 
\[
\delta \colon \mathcal{E}^{\Gamma} (Y,Y',Z)_{< \mathcal{U}}  \longrightarrow
{E}^{\Gamma} (Y,Y',Z),
\]
induced by inclusions $\G_{\Gamma} (Y_S, Y'_S,Z)_{<U} \to \G_{\Gamma} (Y_S, Y'_S,Z)$,
is a weak equivalence.
\end{ThmRefName}

\begin{proof}
	Apply Theorem \refT{Exc2} inductively to the maximal sets in $\mathcal{U}$.
\end{proof}

\SecRef{Other properties of equivariant fibred $G$-theory}{SGP}

\SSecRef{Fibred assembly map}{COP}

The usual notion of metric assumes only finite values.  We will require a \textit{generalized metric} on a set $X$. It is a function
$d \colon X \times X \to [0,\infty) \cup \{\infty\}$
which is reflexive, symmetric, and satisfies the triangle inequality
in the obvious way.
The generalized metric space is \textit{proper} if it is a countable disjoint union of
metric spaces $X_i$ on each of which the generalized metric $d$ is finite,
and all closed metric balls in $X$ are compact.
The metric topology on a generalized metric space is
defined as usual.  
 
The basic fibred assembly map
\[
A (X, Y) \colon \hlf (X; \Gnc (Y))
\longrightarrow \GncX (Y),
\]
for a proper generalized metric spaces $X$ and $Y$, sends the locally finite homology of $X$ with coefficients in the spectrum $\Gnc (Y)$ to the nonconnective fibred $G$-theory $\GncX (Y)$ defined in section \refS{KICVB}.

The locally finite homology $\hlf (X; \mathcal{S})$ we use was introduced in \cite[Definition II.5]{gC:95} for any coefficient spectrum $\mathcal{S}$.
Let $^b S_k X$ be the collection of all locally finite families $\mathcal{F}$
of singular $k$-simplices in $X$ which are uniformly bounded, in the sense that each family possesses a number $N$ such that the diameter of the image $\im (\sigma)$ is bounded from above by $N$ for all simplices $\sigma \in \mathcal{F}$.
For any spectrum $\mathcal{S}$, the theory $\bhlf (X; \mathcal{S})$ is the realization of the simplicial spectrum
\[
  k \mapsto \hocolim{C \in \, ^b \! S_k X} \hlf (C, \mathcal{S}).
\]
There is an equivalence of spectra
$\bhlf (X; \mathcal{S})
\rightarrow
\hlf (X; \mathcal{S})$,
for any proper generalized metric space $X$, from \cite[Corollary II.21]{gC:95}.

A similar theory $\Jh (X,\mathcal{A})$ is obtained as the realization of the
simplicial spectrum
\[
  k \mapsto \hocolim{C \in ^b \! S_k X} \Knc (C,\mathcal{A})
\]
by viewing $C$
as a discrete metric space and using the notation $\Knc (C,\mathcal{A})$ for the nonconnective delooping of the $K$-theory of $\mathcal{B} (C,\mathcal{A})$ from Definition \refD{OPQWS}.
Using the coefficients $\mathcal{A} = \B_C (Y)$, we obtain $\Jh (X,\mathcal{A})$ which we denote $\Jh (X,Y)$.
The proof of \cite[Corollary III.14]{gC:95} gives a
weak homotopy equivalence
\[
  \eta \colon \hlf (X; \Gnc (Y))
  \longrightarrow \Jh (X,Y)
\]
of functors from proper locally compact metric spaces and coarse maps to spectra.

We next define a natural
transformation
\[
  \ell \colon \Jh (X,Y)
  \longrightarrow \GncX (Y).
\]
In the case $Y$ is a point and the coefficients are finitely
generated free $R$-modules, this kind of transformation is defined as part of the proof of Proposition III.20 of \cite{gC:95}.
The definition is entirely in terms of maps between
singular simplices in $X$, so the construction can be generalized to give $\ell$ as above.
For convenience of the reader, we present the necessary details.

Let us first note that controlled algebra can be used to build equivalent bounded $K$-theory spectra using the symmetric monoidal category approach which we will find useful in the rest of the paper. For the details we refer to section 6 of \cite{gC:05}.

Let $\mathcal{D}$ be any collection of singular $n$-simplices of $X$ and $\zeta$ be any point of the standard $n$-simplex.
Define a function $\vartheta_{\zeta} \colon \mathcal{D} \to X$ by $\vartheta_{\zeta} (\sigma) = \sigma (\zeta)$.
Since $\mathcal{D}$ is viewed as a discrete metric space, if $\mathcal{D}$ is locally finite then $\vartheta_{\zeta}$ is coarse, so we
have the induced functor $\mathcal{B} (\mathcal{D}, \mathcal{A}) \to \mathcal{B} (X, \mathcal{A})$ given by
\[
\bigoplus_{d \in \mathcal{D}} F_d \longrightarrow \bigoplus_{x \in X}  \bigoplus_{\vartheta_{\zeta}(d) = x} F_d
\]
which is the identity for each $d \in \mathcal{D}$.
Therefore, there is the induced map of spectra
\[
K(\vartheta_{\zeta}, \mathcal{A}) \colon
K (\mathcal{D}, \mathcal{A}) \longrightarrow K (X, \mathcal{A}).
\]
Suppose further that $\mathcal{D} \in {^b} \!  S_k X$ and that $N$ is a bound required to exist for $\mathcal{D}$ in ${^b} \!  S_k X$.
If $\zeta$ and $\theta$ are both points in the standard $n$-simplex, we have a symmetric monoidal natural transformation
$N_{\zeta}^{\theta} \colon K(\vartheta_{\zeta}, \mathcal{A}) \to K(\vartheta_{\theta}, \mathcal{A})$ induced from the functors which are identities on objects in the cocompletion of $\mathcal{A}$.  Both of those identity morphisms are isomorphisms in
$\mathcal{B} (X, \mathcal{A})$ because they and their inverses are bounded by $N$.

Recall that the standard $n$-simplex can be viewed as the nerve of the ordered set $\underline{n} = \{0, 1, \ldots, n \}$, with the natural order, viewed as a category.
Let $\mathcal{D} \in {^b} \!  S_n X$.  We define a functor
$l(\mathcal{D},n) \colon i\mathcal{B}(\mathcal{D}, \mathcal{A}) \times \underline{n} \rightarrow
i\mathcal{B}(X, \mathcal{A})$
as follows.
On objects, $(l(\mathcal{D},n) F)_x = \bigoplus_{\vartheta(i)=x} F_d$, where $i$ denotes the vertex of $\Delta^n = N . \, \underline{n}$ corresponding to $i$.
On morphisms, $l(\mathcal{D},n)$ is defined by the requirement that the restriction to the subcategory $i\mathcal{B}(\mathcal{D}, \mathcal{A}) \times j$ is the functor induced by $\theta_j$, and that $(\id \times (i \le j)) (F)$ is sent to $N_i^j (F)$.
This is compatible with the inclusion of elements in $^b \!  S_n X$, so we obtain a functor
\[
\colim{\mathcal{D} \in {^b} \!  S_n X} i\mathcal{B}(\mathcal{D}, \mathcal{A}) \times \underline{n}
\longrightarrow
i\mathcal{B}(X, \mathcal{A}),
\]
and therefore a map
\[
\hocolim{\mathcal{D} \in {^b} \!  S_n X} N . \, i\mathcal{B}(\mathcal{D}, \mathcal{A}) \times \Delta^n
\longrightarrow
N . \, i\mathcal{B}(X, \mathcal{A}).
\]
If $\mathcal{M}$ is a symmetric monoidal category, let the $t$-th space in $\Spt (\mathcal{M})$ be denoted by $\Spt_t (\mathcal{M})$,
and let $\sigma_t \colon S^1 \wedge \Spt_t (\mathcal{M}) \to \Spt_{t+1} (\mathcal{M})$ be the structure map for $\Spt (\mathcal{M})$.
The fact that the natural transformations $N_i^j$ are symmetric monoidal shows in particular that we obtain maps
\[
\Lambda_t \colon
\hocolim{\mathcal{D} \in {^b} \!  S_n X} \Spt_t (i\mathcal{B}(\mathcal{D}, \mathcal{A})) \times \Delta^n
\longrightarrow
\Spt_t (i\mathcal{B}(X, \mathcal{A}))
\]
so that the diagrams
\[
\xymatrix{
  \hocolim{\mathcal{D} \in {^b}\!S_n X} (S^1 \wedge \Spt_{t} (i\mathcal{B}(\mathcal{D}, \mathcal{A}))) \times \Delta^n
  \ar[r] \ar[d]_-{\sigma_t \times \id}
  &S^1 \wedge \Spt_t (i\mathcal{B}(X, \mathcal{A})) \ar[d]^-{\sigma_t}\\
  \hocolim{\mathcal{D} \in {^b} \!  S_n X} \Spt_{t+1} (i\mathcal{B}(\mathcal{D}, \mathcal{A})) \times \Delta^n
   \ar[r]^-{\Lambda_{t+1}}
  &\Spt_{t+1} (i\mathcal{B}(X, \mathcal{A}))
}
\]
commute.
Further, for each $t$ we obtain a map
\[
\left| \, \underline{k} \mapsto \hocolim{\mathcal{D} \in {^b} \!  S_n X} \Spt_{t} (i\mathcal{B}(\mathcal{D}, \mathcal{A})) \, \right|
\longrightarrow
\Spt_t (i\mathcal{B}(X, \mathcal{A})) 
\]
respecting the structure maps in $\Spt_t$.
This gives a map
$\ell \colon {^c} \!\! \Jh (X; \mathcal{A})
  \rightarrow K (X, \mathcal{A})$
where ${^c} \!\! \Jh (X; \mathcal{A})$ stands for the realization of the
simplicial spectrum
\[
  k \mapsto \hocolim{C \in {^b} \! S_k X} K (C, \mathcal{A}).
\]
Since $\ell$ is natural in $X$ and is compatible with delooping,
it generalizes to the homotopy natural transformation
$\ell_K \colon \Jh (X; \mathcal{A})
  \rightarrow \Knc (X, \mathcal{A})$.
Composing this with the Cartan natural transformation $\Knc (X, \mathcal{A}) \to \GncX (Y)$ gives 
\[
\ell \colon \Jh (X,Y)
  \longrightarrow \GncX (Y).
  \]

\begin{DefRefName}{GYPPL}{Fibred assembly map in $G$-theory}
The homotopy natural transformation
\[
  A (X,Y) \colon
  \hlf (X; \Gnc (Y))
  \longrightarrow \GncX (Y)
\]
is the composition of $\eta$ and $\ell$.
\end{DefRefName}

\SSecRef{Fibrewise trivialization}{COPSS}

In this section we want to justify the claim from the introduction that in the new equivariant theory we have built there are fibrewise trivializations.  First, we will state the desired fact precisely.

Recall the proper metric space $Y_S$ described in section \refS{TTITS}.  We assume that the trivial action $s_0$ is in $S$ and use the notation $Y_0$ for the space $Y$ with the trivial action.

\begin{ThmRef}{KTDHHJDD}
	The equivariant inclusion of metric spaces $Y_0 \to Y_S$ induces an equivalence $\widetilde{G}^{\Gamma} (Y_0)^{\Gamma} \to \widetilde{G}^{\Gamma} (Y_S)^{\Gamma}$.  Therefore, there is an equivalence $\widetilde{G}^{\Gamma} (Y_0) \to E^{\Gamma} (Y)$. 
\end{ThmRef}

We start with several facts about filtrations on modules.
Let $\Phi_d \colon \mathcal{P} (X) \to \mathcal{P} (X)$ denote the functor that assigns to a subset of $X$ its $d$-neighborhood in $X$.  We can think of an object of $\G (X)$ as a pair $(F,\theta)$, where $\theta$ has a grading with the properties spelled out in Definition \refD{LScov22}.  Given two $X$-filtrations $\theta$ and $\eta$, we say $\theta$ is contained in $\eta$ if $\theta (S) \subset \eta (S)$ for all $S \subset X$ and write $\theta \le \eta$. We say two $X$-filtrations $\theta$ and $\eta$ are similar if there is a number $d$ so that
$\theta \le \eta \circ \Phi_d$ and $\eta \le \theta \circ \Phi_d$.  

\begin{LemRef}{NVDXC}
If $\theta$ and $\eta$ are similar then the objects $(F,\theta)$ and $(F,\eta)$ are isomorphic in $\G (X)$.
\end{LemRef}

\begin{proof}
The conditions ensure that the identity homomorphism is boundedly controlled in both directions.	
\end{proof}

Let $f\colon X \to Y$ be a coarse map of proper metric spaces, so that for every $d \ge 0$ there is $L(d) \ge 0$ so that $d(x,y) \le d$ implies $d(f(x),f(y)) \le  L(d)$. Suppose further that $f$ is proper. Given an $X$-filtration $\theta$ on an $R$-module $F$, we define $f_{\ast} (\theta)$ to be the $Y$-filtration on $F$ given by $f_{\ast} (\theta)(U) = \theta (f^{-1} (U))$. Similarly, given a $Y$-filtration $\theta$ on $F$, we define $f^{\ast} (\theta)$ to be the $X$-filtration on $F$ given by $f^{\ast} (\theta)(U) = \theta (f (U))$.

Recall the definition of $Y_S$ in \refD{MPEW}.  
Let $i \colon Y \to Y_S$ be the inclusion $y \to  (y,s_0)$, an isometric embedding, and let $\pi \colon Y_S \to Y$ denote the projection, a distance non-increasing map.

\begin{LemRef}{NVDXCC}
Let $F$ be any $R$-module.  Then any $Y_S$-filtration on $F$ is similar to one of the form $i_{\ast} \theta$, where $\theta$ is a $Y$-filtration on $F$.
\end{LemRef}

\begin{proof}
Let $ex \colon \mathcal{P}(Y_S) \to \mathcal{P}(Y_S)$ be defined by $ex(U) = \pi(U) \times S$. It is clear from the definition that $U \subset ex(U)$. It is also readily checked that $ex(U) \subset \Phi_1 (U)$, which shows that any $Y_S$-filtration $\theta$ on an $R$-module $F$ is similar to the $Y_S$-filtration $\theta \circ ex$. Let $\overline{\theta}$ denote the $Y$-filtration on $F$ given by $\overline{\theta}(U) = \theta (U \times S)$. 
Then it is clear that $\theta \circ ex = \pi^{\ast} (\overline{\theta})$. 
It therefore suffices to show that for any
$Y$-filtration $\eta$ on $F$, 
we have that $i_{\ast} \eta$ and $\pi^{\ast} \eta$ are similar. 
But it is clear that $\pi^{\ast} \eta \le i_{\ast} \eta \circ \Phi_1$, which gives the result.	
\end{proof}

We also have the following useful fact.

\begin{LemRef}{NVDXCCC}
Suppose that we are given two $X$-filtrations $\theta$ and $\eta$ on an $R$-module $F$, and that $f \colon F \to F$ is bounded as a morphism from $(F,\theta)$ to $(F,\eta)$. Suppose further that $\theta'$ and $\eta'$ are also $X$-filtrations, and that $\theta'$ and $\eta'$ are similar to $\theta$ and $\eta$, respectively. Then $f \colon F \to F$ is bounded as a morphism from from $(F,\theta')$ to $(F,\eta')$.
\end{LemRef}

Suppose a metric space $Y$ has an action by a discrete group $\Gamma$ through coarse equivalences.  
Recall that we say \textit{the action is bounded} if for each
$\gamma \in \Gamma$, 
there is 
$b(\gamma) \ge 0$ 
so that $d(y, \gamma y) \le b(\gamma)$
for all 
$y \in Y$. 
The following is an elementary observation.

\begin{LemRef}{NVDXCCCC}
	Suppose $Y$ is a 
	proper metric space equipped with a bounded $\Gamma$-action by coarse maps.
	Then, given any $Y$-filtration $\theta$ on an $R$-module $F$, and any $\gamma \in \Gamma$
	we have that $\theta$ and $\gamma_{\ast} \theta$ are similar.
\end{LemRef}

We have the following equivalent interpretation of the category $\G^{h\Gamma} (Y)$ introduced in Definition \refD{MPEW}.  An object of $\G^{h\Gamma} (Y)$ is given by data $(F,\theta, \{ f_{\gamma,\gamma'} \} _{\gamma, \gamma' \in \Gamma})$ where
\begin{enumerate}
	\item $F$ is an $R$-module,
	\item $\theta$ is a $Y$-filtration on $F$,
	\item $f_{\gamma, \gamma'}$ is an automorphism of $F$,
	\item $f_{\gamma,\gamma} = \id_F$ and $f_{\gamma, \gamma'} \circ f_{\gamma', \gamma''} = f_{\gamma, \gamma''}$,
	\item $f_{\gamma, \gamma'}$ is bounded when regarded as a homomorphism $(F,\gamma'_{\ast} \theta) \to (F,\gamma_{\ast} \theta)$.
\end{enumerate}

Lemmas \refL{NVDXCCC} and \refL{NVDXCCCC} give that condition 5 on $f_{\gamma, \gamma'}$ is equivalent to  $f_{\gamma, \gamma'}$ being bounded as a homomorphism from $(F,\theta)$ to $(F,\theta)$.

Now we are ready to prove Theorem \refT{KTDHHJDD}.

\begin{proof}
First, observe it suffices to verify that the inclusion induces an equivalence of categories $\iG^{h\Gamma} (Y_0) \to \iG^{h\Gamma} (Y_S)$.  The equivalence then clearly extends to categories of diagrams of objects in $\G^{h\Gamma} (Y_S)$, and Waldhausen's $S.$-construction used to produce the spectra gives simplicial spaces which in every level are the nerves of categories of isomorphisms of diagrams of cofibrations of objects in $\G^{h\Gamma} (Y_S)$. 

The inclusion exhibits $\iG^{h\Gamma} (Y_0)$ as a full subcategory of $\iG^{h\Gamma} (Y_S)$, and it follows that it's enough to prove that every object of $\iG^{h\Gamma} (Y_S)$ is isomorphic to an object of $\iG^{h\Gamma} (Y_0)$.
An object of $\iG^{h\Gamma} (Y_0)$ is given by data 
$(F,\theta, \{ f_{\gamma,\gamma'} \} _{\gamma, \gamma' \in \Gamma})$,
where $\theta$ is an $Y_0$-filtration on $F$, and where $f_0$ is an automorphism of $F$ which is bounded as a homomorphism from $(F,\theta)$ to $(F,\theta)$. 
Note that the transformations by $\gamma$'s do
not occur in this situation because the 
action of $\Gamma$ on $Y_0$ is trivial. The inclusion functor $\iG^{h\Gamma} (Y_0) \hookrightarrow \iG^{h\Gamma} (Y_S)$ is given by 
\[
(F,\theta, \{ f_{\gamma,\gamma'} \} _{\gamma, \gamma' \in \Gamma}) \to (F, i_{\ast}\theta, \{ f_{\gamma,\gamma'} \} _{\gamma, \gamma' \in \Gamma}),
\]
so an object $(F,\theta, \{ f_{\gamma,\gamma'} \} _{\gamma, \gamma' \in \Gamma})$ is in the subcategory $\iG^{h\Gamma} (Y_0)$ if and only if $\theta$ is of the form $i_{\ast} \eta$ for some $Y_0$-filtration $\eta$.

Next, we observe that if $(F,\theta, \{ f_{\gamma,\gamma'} \} _{\gamma, \gamma' \in \Gamma})$ is an object of $ \iG^{h\Gamma} (Y_S)$, and if  $\theta'$ is an $Y_S$-filtration on $F$ which is similar to $\theta$, then (a) $(F,\theta', \{ f_{\gamma,\gamma'} \} _{\gamma, \gamma' \in \Gamma})$ is also an object of $(F, i_{\ast}\theta, \{ f_{\gamma,\gamma'} \} _{\gamma, \gamma' \in \Gamma})$, and (b) $(F, \theta, \{ f_{\gamma,\gamma'} \} _{\gamma, \gamma' \in \Gamma})$ is isomorphic to $(F, \theta', \{ f_{\gamma,\gamma'} \} _{\gamma, \gamma' \in \Gamma})$. 
But we have already observed in Lemma \refL{NVDXCC} that every $Y_S$-filtration on $F$ is equivalent to one of the form $i_{\ast} \eta$, for some $Y_0$-filtration $\eta$ on $F$, proving the result.	
\end{proof}

\end{document}